\documentclass[12pt,a4paper]{article}

\usepackage[francais]{babel}
\usepackage{amsfonts}
\newcommand {\CO}{{\cal O}}
\newcommand{\C}{\mathbb C}
\newcommand{\D}{\mathbb D}
\newcommand{\G}{\mathbb G}
\newcommand{\Hh}{\mathbb H}
\newcommand{\N}{\mathbb N}
\newcommand{\R}{\mathbb R}
\newcommand{\Sp}{\mathbb S}
\newcommand{\U}{\mathbb U}
\newcommand{\W}{\mathbb W}
\newcommand{\pv}{{\it D\'emonstration : }}
\newcommand{\vp}{\hfill $ \dagger $}

\newtheorem{theor}{Th\'eor\`eme}[section]
\newtheorem{conj}[theor]{Conjecture}
\newtheorem{deft}[theor]{D\'efinition}
\newtheorem{ppt}[theor]{Propri\'et\'e}

\newtheorem{lem}[theor]{Lemme}
\newtheorem{prop}[theor]{Proposition}
\newtheorem{cor}{Corollaire}[theor]
\newtheorem{rem}[cor]{Remarque}

\newif\ifpdf
\ifx\pdfoutput\undefined
\pdffalse                   % we are not running PDFLaTeX
\else
\pdfoutput=1                % we are running PDFLaTeX
\pdftrue
\fi

\ifpdf
    \pdfoutline goto name {titre} {Titre}
    \pdfcatalog{/PageMode /UseOutlines}
   \usepackage[pdftex]{color}
   \usepackage[pdftex,colorlinks,a4paper,hyperindex]{hyperref}
\else
   \usepackage{hyperref}
\fi

\begin{document}

    \ifpdf \pdfdest name {titre} XYZ \fi
    
    \author{Beno\^{\i}t RIVET \thanks{
            L'auteur remercie les soci\'et\'es Altran
	    Technologies, France Telecom Mobiles et FirstMark
	    Communications France pour le cadre agr\'eable
	    qu'elles lui ont fourni au cours de ce travail.    
	    Adresse \'electronique : 
	    \texttt{benoit.rivet@online.fr}}}

    \title{Dif\-f\'eo\-mor\-phis\-mes harmoniques du plan 
           hyperbolique}
    
    \maketitle
    
    \begin{abstract}
	Sampson \cite{Sampson-78} et Schoen et Yau \cite{Schoen-Yau-78} ont 
	d\'emontr\'e que tout dif\-f\'eo\-mor\-phis\-me entre surface 
	de Riemann hyperboliques est homotope \`a un dif\-f\'eo\-mor\-phis\-me
	harmonique. Comme l'avait conjectur\'e Schoen \cite{Schoen-93}
	et comme l'avaient partiellement d\'emontr\'e Wan \cite{Wan-92} 
	et Tam et Wan \cite{Tam-Wan-95}, nous d\'emontrons dans cet article
	que ce r\'esultat se g\'en\'eralise aux surfaces non compactes :
	tout ho\-m\'eo\-mor\-phis\-me quasi sym\-m\'e\-tri\-que du cercle 
	s'\'etend en un dif\-f\'eo\-mor\-phis\-me harmonique quasi
	conforme du plan hyperbolique.
	Ce th\'eor\`eme permet de donner une param\'etrisa\-tion de 
	l'espace universel de Teichm\"uller par les diff\'erentielles
	quadratiques holomorphes born\'ees sur le plan hyperbolique.
    \bigskip
    
    \begin{center}
	{\bf Abstract}
    \end{center}
	A classical result of Sampson \cite{Sampson-78} and Schoen and Yau 
	\cite{Schoen-Yau-78} states that every diffeomorphism between 
	compact hyperbolic Riemann surfaces is homotopic to an harmonic
	diffeomorphism.  As conjectured by Schoen \cite{Schoen-93} and 
	partially proved by Wan \cite{Wan-92} and Tam and Wan 
	\cite{Tam-Wan-95}, we prove in this article that this theorem 
	generalizes to the non compact case :
	every quasi symmetric homeomorphism of the circle extends
	to a quasi isometric harmonic diffeomorphism of the hyperbolic plane.
	This enables to parametrize the universal Teichm\"uller space by 
	bounded holomorphic quadratric differentials of the hyperbolic plane
    \end{abstract}
    
    \bigskip

    La dif\-f\'e\-ren\-tiel\-le de Hopf d'une application
    entre deux surfaces de Riemann $ f $ : $ M\rightarrow N $ 
    est la dif\-f\'e\-ren\-tiel\-le quadratique 
    de type $ (2,0) $ :  $ \phi=\varphi(z)dz^2 $
    d\'ecrivant la partie sans trace de $ f^\ast g_{N} $ :
    $$ f^\ast g_{N} = \phi
    + \frac{\parallel df \parallel^{2}}{2} g_{M} + \bar{\phi} $$
    La diff\'erentielle de Hopf est un outil essentiel pour 
    l'\'etude des applications harmoniques entre surfaces de
    Riemann. En effet, lorsque $ f $ est une application harmonique,
    sa dif\-f\'e\-ren\-tiel\-le de Hopf est holomorphe.
    R\'eciproquement, si la diff\'erentielle de 
    Hopf de $ f $ est holomorphe et si $ f $ est injective, $ f $ est 
    une application harmonique.
    
    Le lien entre dif\-f\'eo\-mor\-phis\-mes harmoniques et
    diff\'erentielle de 
    Hopf est enti\`erement r\'esolu pour les surfaces de Riemann 
    hyperboliques compactes. Sampson \cite{Sampson-78} et
    Schoen et Yau \cite{Schoen-Yau-78} ont en effet d\'emontr\'e 
    le th\'eor\`eme d'existence :
    \begin{quote}
	\it Si $ f $ est un dif\-f\'eo\-mor\-phis\-me entre deux surfaces 
	    de Riemann compactes hyperboliques, $ f $ est homotope \`a un
	    unique dif\-f\'eo\-mor\-phis\-me harmonique.
    \end{quote}
    et Sampson \cite{Sampson-78} a d\'emontr\'e le th\'eor\`eme d'unicit\'e :
    \begin{quote}
	\it Si $ M $ et $ N $ sont des surfaces de Riemann hyperboliques 
	compactes,
	si $ f_{1} : M\rightarrow N $ et $ f_{2} : M\rightarrow N $
	sont deux dif\-f\'eo\-mor\-phis\-mes harmoniques qui ont m\^eme
	diff\'erentielle de Hopf, alors $ f_{1}\circ f_{2}^{-1} $ est 
	une isom\'etrie.
    \end{quote}
    Sampson \cite{Sampson-78} en d\'eduisait que 
    l'espace de Teichm\"uller des surfaces compactes s'injecte dans 
    l'espace des dif\-f\'e\-ren\-tiel\-les quadratiques holomorphes.
    Si on fixe
    une structure hyperbolique de r\'ef\'erence $ g_{0} $ sur une 
    surface de Riemann compacte $ M $, on peut associer \`a toute 
    m\'etrique hyperbolique $ g $ sur $ M $ la 
    dif\-f\'e\-ren\-tiel\-le de Hopf de l'unique 
    dif\-f\'eo\-mor\-phis\-me harmonique
    $$ f  :  (M,g_{0})\rightarrow (M,g) $$ 
    homotope \`a l'identit\'e. Les travaux de
    Wolf \cite{Wolf-89} ont permis de v\'erifier que l'application 
    ainsi d\'efinie est bijective : \`a toute diff\'erentielle de 
    Hopf est associ\'ee une unique m\'etrique hyperbolique $ g $
    telle que l'application identit\'e
    $ id $ : $ (M,g_{0})\rightarrow (M,g) $ 
    est harmonique. Les travaux de Sampson et Wolf permettent de 
    conclure :    
    \begin{quote}
	\it Soit $ M $ une surface de Riemann compacte.
	Si $ \varphi $ est une diff\'eren\-tiel\-le quadratique holomorphe sur 
	$ M $ $ \varphi $,
	il existe une unique m\'e\-tri\-que hyperbolique $ g_{\varphi} $
	telle que l'application identit\'e
	$$ id :(M,g_{0})\rightarrow (M,g_{\varphi}) $$
	est harmonique.
	L'application qui \`a $ \varphi $ associe la classe de Teich\-m\"uller 
	de la m\'etrique $ g_{\varphi} $ d\'efinit 
	une bijection entre l'espace des diff\'erentielles quadratiques 
	holomorphes et l'espace de Teichm\"ul\-ler de M.
    \end{quote}
    
    Ces r\'esultats ont \'et\'e partiellement
    g\'en\'eralis\'es \`a l'espace universel de 
    Teichm\"uller par Wan \cite{Wan-92} et Tam et Wan \cite{Tam-Wan-95}.
    Wan \cite{Wan-92} a d\'emontr\'e l'existence d'une application de 
    l'espace des diff\'erentielles quadratiques holomorphes born\'ees 
    dans l'espace de Teichm\"uller :
    \begin{quote}
	\it A toute dif\-f\'e\-ren\-tiel\-le quadratique holomorphe
	$ \phi $ born\'ee sur $ \Hh^2 $ est associ\'ee un
	dif\-f\'eo\-mor\-phis\-me du plan hyperbolique 
	$ f $ de dif\-f\'e\-ren\-tiel\-le de Hopf $ \phi $.
	$ f $ est unique \`a isom\'etrie pr\`es et est une 
	quasi isom\'etrie.
    \end{quote}
    et Tam et Wan \cite{Tam-Wan-95} ont d\'emontr\'e une r\'eciproque 
    partielle :   
    \begin{quote}
	\it L'application $ QD_{b}(\Hh^2) \rightarrow Teich(\Hh^2) $ qui \`a 
	une diff\'erentielle quadratique born\'ee $ \phi $ associe la classe de 
	Teichm\"uller des dif\-f\'eo\-mor\-phis\-mes harmoniques quasi conforme
	de diff\'erentielle de Hopf $ \phi $ est d'image ouverte.
    \end{quote}

    Le probl\`eme de la surjectivit\'e de l'application d\'efinie par Wan 
    est \'equi\-va\-lent \`a la conjecture formul\'ee par Schoen 
    \cite{Schoen-93} :
    \begin{conj}[Schoen]
	\it Tout ho\-m\'eo\-mor\-phis\-me quasi sym\-m\'e\-tri\-que du 
	cer\-cle s'\'etend en un 
	dif\-f\'eo\-mor\-phis\-me harmonique quasi conforme
	du plan hyperbolique.
    \end{conj}
    
    Nous d\'emontrons dans cet article la conjecture de Schoen. Le 
    plan de l'article est le suivant :
    
    \tableofcontents

    Dans la premi\`ere partie, nous rappelons quelques 
    r\'esultats \'el\'ementaires sur les surfaces de Riemann,
    l'espace de Teichm\"uller et les
    propri\'et\'es de stabilit\'e et de compacit\'e des
    applications harmoniques.

    Dans la seconde partie, nous d\'efinissons la notion de 
    m\'etrique harmonique et nous donnons une exposition 
    d\'etaill\'ee de la th\'eorie 
    de Wan et de la r\'esolution de l'\'equation des 
    m\'etriques harmoniques :
    $$ \frac{1}{2}\Delta h = e^h-|\phi|^2 e^{-h}-1 $$
    Nous donnons en particulier des exemples de 
    dif\-f\'e\-ren\-tiel\-les de 
    Hopf (non born\'ees) sur $ \Hh^{2} $ qui v\'erifient la 
    propri\'et\'e de non-compl\'etude :
    \begin{ppt}
	Si $ \phi $ est la diff\'erentielle de Hopf d'une
	immersion harmonique $ f $, 
	$ f $ n'est pas surjective.
    \end{ppt}
    
    Ces exemples ont d\'ej\`a \'et\'e remarqu\'es par Shi, Tam et Wan 
    \cite{Shi-Tam-Wan-99}. L'analyse d\'etaill\'ee de ces exemples 
    permet de mettre en \'evidence le r\^ole jou\'e par 
    la g\'eom\'etrie de la diff\'erentielle de Hopf et des feuilletages 
    associ\'es pour \'etudier la
    surjectivit\'e de $ f $.
    \`A l'exception d'un principe du maximum pour les facteurs 
    de distortion quasi conforme des matriques harmoniques et de 
    l'\'etude d\'etaill\'ee des exemples de Shi, Tam et Wan 
    \cite{Shi-Tam-Wan-99}, les r\'esultats expos\'es dans cette 
    partie sont ceux de Wan \cite{Wan-92}, 
    expos\'es suivant un point de vue l\'eg\`erement diff\'erent 
    (d\'ej\`a explicit\'e dans ma th\`ese \cite{Rivet-99}).
    
    Nous d\'emontrons dans la troisi\`eme partie le 
    th\'eor\`eme principal de notre \'etude :
    \begin{theor}
	Tout ho\-m\'eo\-mor\-phis\-me quasi sym\-m\'e\-tri\-que du 
	cercle s'\'etend en un 
	dif\-f\'eo\-mor\-phis\-me harmoni\-que quasi conforme
	du plan hyperbolique.
    \end{theor}

    \section{Applications harmoniques des surfaces de Riemann}
    
    \subsection{Surfaces de Riemann}
    
    Nous rappelons dans ce paragraphe quelques propri\'et\'es 
    \'el\'ementaires des surfaces de Riemann.

    Une structure complexe $ J $ sur une vari\'et\'e
    $ M $ de dimension 2 est une section de 
    $ T^\ast M \otimes TM $ v\'erifiant $ J\circ J = -id $.
    Si $ J $ est une structure complexe sur $ M $,
    $ (M,J) $ est une surface de Riemann.
    Si $ M $ et $ N $ sont deux surfaces de 
    Riemann, une application $ f $~: $ M\rightarrow N $
    est holomorphe si et seulement si $ f $ commute avec les 
    structures complexes de $ M $ et $ N $~: 
    $ J_N\circ df = df\circ J_M $.
    
    Une m\'etrique $ g $ sur $ M $ est une
    m\'etrique conforme
    si la conjugaison complexe $ J $ laisse $ g $ invariant : $ g\circ J = g $.
    
    L'\'etude des surfaces de Riemann est bas\'ee sur le 
    th\'eor\`eme d'uniformisation de Riemann~:

    \begin{theor}
	Soit $ (M,J) $ une surface de Riemann.
	Le rev\^etement universel de $ (M,J) $ est biholomorphe :
	\begin{itemize}
	    \item soit \`a la sph\`ere de Riemann $ \Sp^2 $,
	    \item soit au plan $ \C $,
	    \item soit au demi plan $ \U=\{x+iy|y>0\} $.
	\end{itemize}
	$ M $ est elliptique lorsque son rev\^etement universel est
	biholomorphe \`a $ \Sp^{2} $, parabolique lorsqu'il est biholomorphe
	\`a $ \C $ et hyperbolique lorsqu'il est biholomorphe \`a $ \U $. 
    \end{theor}

    Lorsque $ M $ est hyperbolique, la donn\'ee d'une structure 
    conforme est \'equi\-va\-len\-te \`a la donn\'ee d'une orientation
    et d'une m\'etrique hyperbolique.
    Le groupe des biholomorphismes de $ \U $ s'identifie en effet \`a 
    $ PSL_2(\R) $, agissant par homographie.
    Il laisse invariant une unique famille 
    de m\'etriques conformes~:
    $$ g_t=t\frac{dx^2+dy^2}{y^2} $$
    de courbure $ K_t=-t^{-1} $.
    Si $ (M,J) $ est une surface de Riemann hyperbolique,
    $ M $ est donc canonniquement muni d'une m\'etrique 
    conforme \`a courbure $ -1 $.
    C'est la m\'etrique de Poincar\'e de $ (M,J) $.
    
    Si $ f $ : $ (M,g_{M})\rightarrow(N,g_{N}) $ est une application 
    $ C^1 $ entre deux surfaces de Riemann munies de m\'etriques 
    conformes, on notera $ (\lambda_{1}(f))^2 $ et $ (\lambda_{2}(f))^2 $
    les valeurs propres de $ f^\ast g_{N} $ par rapport \`a
    $ g_{M} $. $ \lambda_{1} $ et $ \lambda_{2} $ sont les 
    coefficients de dilatation de $ f $ et sont rang\'es par ordre 
    d\'ecroissant : $ \lambda_{1} \geq \lambda_{2} \geq 0 $.
    $ f $ est $ K $ quasi isom\'etrique si ses coefficients de 
    dilatation sont born\'es :
    $$ K \geq \lambda_{1} \geq \lambda_{2} \geq \frac{1}{K} $$
    $ f $ est $ K $ quasi conforme si le rapport de quasi 
    conformit\'e $ \frac{\lambda_{1}}{\lambda_{2}} $ est born\'e :
    $$ \frac{\lambda_{1}}{\lambda_{2}} \leq K $$
    Le rapport de quasi conformit\'e 
    $ K(f)=\sup_{M}\frac{\lambda_1(f)}{\lambda_2(f)} $
    ne d\'epend pas du choix des
    m\'etriques con\-for\-mes sur $ M $ et $ N $, mais seulement des 
    structures complexes $ J_{M} $ et $ J_{N} $ de $ M $ et $ N $.
    Deux structures complexes $ J_{1} $ et $ J_{2} $ d\'efinies sur 
    une m\^eme surface de Riemann $ M $ 
    sont quasi conformes si l'application identit\'e de $ (M,J_{1}) $ 
    vers $ (M,J_{2}) $ est quasi conforme.
    L'espace des structures conformes $ J $
    quasi conformes \`a une structure de r\'ef\'erence $ J_{M} $ est 
    naturellement muni de la distance de Teichm\"uller :
    $$ d(J_{1},J_{2}) = \log(K(id)) $$
    o\`u $ K(id) $ est le rapport de quasi conformit\'e de 
    l'application identit\'e entre $ (M,J_{1}) $ et $ (M,J_{2}) $.
    De m\^eme, l'espace des m\'etriques hyperboliques est 
    naturellement muni de la distance de Teichm\"uller.
    
    L'espace universel de Teichm\"uller est l'espace de Teich\-m\"uller 
    du plan hyperbolique $ \Hh^2 $, muni de la structure conforme de 
    r\'ef\'erence $ J_{\Hh^2} $ et de la m\'etrique de Poincar\'e 
    $ g_{\Hh^2} $. On peut le d\'efinir comme :
    \begin{itemize}
	\item l'espace des modules des structures complexes $ J $
	      sur $ \Hh^2 $, muni de la distance de Teichm\"uller,
	      modulo l'action du groupe des hom\'eomorphismes quasi conforme
	      \`a distance born\'ee de l'identit\'e,
	\item l'espace des modules des m\'etriques hyperboliques, muni de la 
	      distance de Teichm\"uller,
	      modulo l'action du groupe des hom\'eomorphismes quasi conformes
	      \`a distance born\'ee de l'identit\'e,
	\item l'espace des hom\'eomorphismes quasi conforme de $ \Hh^2 $,
	      muni de la distance de Teichm\"uller,
	      modulo l'action du groupe des hom\'eomorphismes quasi conformes
	      \`a distance born\'ee de l'identit\'e.
    \end{itemize}
    
    Nous avons d\'ej\`a explicit\'e le lien entre structure complexe 
    et m\'etrique hyperbolique : \`a $ J $ est associ\'ee la m\'etrique 
    de Poincar\'e et r\'eciproquement, $ \Hh^{2} $ \'etant d\'ej\`a 
    orient\'e, \`a toute m\'etrique hyperbolique est associ\'ee une 
    structure complexe.

    D'apr\`es le th\'eor\`eme 
    de Beltrami, si $ g $ est une m\'etrique quasi conforme,
    $ (\Hh^2,g) $ est isom\'etrique \`a $ (\Hh^2,g_{\Hh^2}) $. 
    On associe donc naturellement \`a $ g $ les isom\'etries
    $$ f : (\Hh^2,g)\rightarrow (\Hh^2,g_{\Hh^2}) $$
    qui sont des hom\'eomorphismes quasi conformes de
    $ \Hh^2 $ muni de la m\'etrique de Poincar\'e $ g_{\Hh^2} $ dans lui 
    m\^eme. R\'eciproquement, \`a un hom\'eomorphisme quasi conforme de $ \Hh^2 $
    $ f $ est associ\'e la m\'etrique quasi conforme $ f^\ast g_{\Hh^2} $.
    
    On dispose d'une quatri\`eme caract\'erisation de l'espace 
    universel de Teich\-m\"uller :
    \begin{itemize}
	\item l'espace universel de Teichm\"uller est l'espace des 
	      hom\'eomorphismes quasi symm\'etriques du bord \`a l'infini de
	      $ \Hh^{2} $.
    \end{itemize}
    
    La structure quasi conforme de $ \Hh^{2} $ induit en effet une 
    structure quasi symm\'etrique sur le bord \`a l'infini, que nous 
    allons expliciter rapidement. Dans le 
    mod\'ele du demi plan $ \U = \{x+iy|y>0\} $, le bord \`a l'infini 
    de $ \Hh^{2} $ s'identifie \`a $ \R\cup\{\infty\} $. Le groupe 
    d'isom\'etrie de $ \Hh^{2} $, $ PSL_{2}(\R) $, agit par homographie 
    sur $ \R\cup\{\infty\} $. Cette action est transitive sur les 
    tripl\'es de points. Si on se fixe quatre points $ (a,b,c,d) \in 
    \partial_{\infty} \U $, il existe une unique homographie $ g $ 
    telle que $ g(a,b,c) = (0,1,\infty) $ et on d\'efinit le 
    birapport :
    $ [a,b,c,d] := g(d) $. Lorsque $ (a,b,c,d) \in \R^{4} $, on a la formule :
    $$ [a,b,c,d] = \frac{d-a}{d-c}\times\frac{b-c}{b-a} $$
    Un hom\'eomorphisme de $ \partial_{\infty}\U $ est quasi symm\'etrique 
    s'il laisse presque invariants les birapports, c'est \`a dire s'il 
    existe $ K $ tel que, pour tout quadrupl\'e $ (a,b,c,d) $, on a :
    $$ \frac{1}{K} \leq \frac{[f(a),f(b),f(c),f(d)]}{[a,b,c,d]} \leq K $$
    
    Si $ f $ est un hom\'eomorphisme quasi conforme de $ \Hh^{2} $,
    $ f $ s'\'etend en $ \partial f $,
    un hom\'eomorphisme quasi symm\'etrique du bord
    \`a l'infini de $ \Hh^2 $.
    De plus, si $ \dot{f} $ est \`a distance born\'ee de $ f $, son 
    extension $ \partial \dot{f} $ co\"{\i}ncide avec $ \partial f $. 
    A toute classe de Teichm\"uller d'hom\'eomorphismes de $ \Hh^{2} $ 
    est donc associ\'e l'extension au bord \`a l'infini de l'un de ses 
    repr\'esentants.
    R\'eciproquement, tout hom\'eomorphisme quasi symm\'etrique de
    $ \partial_{\infty} \Hh^{2} $ s'\'e\-tend en un diff\'eomorphisme quasi 
    conforme de $ \Hh^{2} $ : l'extension de Earle et Eells ou 
    l'extension de Douady Earle \cite{Douady-Earle} conviennent.
    
    \subsection{Applications harmoniques}

    Si $ M $ une vari\'et\'e de dimension 2,
    le choix d'une m\'etrique riemannienne et d'une orientation sur 
    $ M $ d\'eterminent une structure complexe sur $ M $,
    telle que dans une carte holomorphe, la m\'etrique $ g $ 
    s'\'ecrive $ g=M^2|dz|^2 $.
    R\'eciproquement, la donn\'ee d'une structure complexe 
    induit une structure conforme $ [g] $ et une orientation.
    La structure conforme associ\'ee \`a $ J $ est
    celle des m\'etrique localement de la forme $ M^2|dz|^2 $
    dans une carte holomorphe. Lorsque $ M $ est compacte, une 
    application harmonique $ f $ est un point critique pour l'\'energie :
    $$ E(f) = \int_{M} |df|^{2}dx $$
    Localement, les applications harmoniques sont caract\'eris\'ee par 
    l'annulation du champ de tension $ \tau(f) = tr(\nabla df) $.
    
    \begin{deft}
	Soient $ (M,g_{M}) $ et $ (N,g_{N}) $ deux vari\'et\'es riemaniennes.
	Une application $ C^{2} $ $ f $ : $ M\rightarrow N $ est harmonique
	si et seulement si le champ de tension de $ f $ s'annule :
	$$ \tau(f) := tr(\nabla df) = 0 $$
    \end{deft}
    
    Lorsque $ M $ est de dimension $ 2 $,
    l'harmonicit\'e d'une application 
    $ f $~: $ M\rightarrow N $
    est une propri\'et\'e conforme : elle ne d\'epend que du choix
    de la structure complexe $ J $ sur $ M $.
    En effet, la 2-forme de tension~:
    $ \tau(f)dx $
    est ind\'ependante de la m\'etrique conforme choisie sur $ M $.
    
    De plus, la diff\'erentielle de Hopf de $ f $ est holomorphe      
    lorsque $ f $ est harmonique. Lorsque $ N $ est une surface de 
    Riemann, l'holomorphie de la diff\'erentielle de Hopf suffit 
    presque \`a caract\'eriser les applications harmoniques.
    
    Si $ \psi $~: $ M\rightarrow\R $ est une fonction 
    harmonique (non constante) et $ \gamma $~: $ \R\rightarrow N $ est 
    une courbe param\'etr\'ee par l'arc (i.e~: $ |\dot{\gamma}|=1 $),
    la diff\'erentielle de Hopf de $ f=\gamma\circ\psi $ est 
    holomorphe, mais $ f $ n'est harmonique que si $ \gamma $ est une 
    g\'eod\'esique.
    Si on suppose que $ f $ est une immersion, on peut d\'emontrer par 
    contre :
    
    \begin{theor}
	Soient $ M $ et $ N $ deux vari\'et\'es riemanniennes de 
	dimension 2 et $ f $ une immersion de $ M $ dans $ N $.
	$ f $ est harmonique si et seulement si sa diff\'erentielle de 
	Hopf est holomorphe.
    \end{theor}
    
    \subsection{Propri\'et\'es de compacit\'e et de 
                stabilit\'e des applications harmoniques}
    
    L'\'equation des applications harmoniques $ \tau(f)=0 $ est une 
    \'equation elliptique. D'apr\`es la th\'eorie de Schauder (cf Gilbarg 
    et Trudinger \cite{Gilbarg-Trudinger}), les 
    topologies $ C^{1} $ et $ C^\infty $ sont \'equivalente pour les 
    applications harmoniques. Pour les applications harmoniques du 
    plan hyperboliques, on a donc des estim\'ees $ C^\infty $ :
    
    \begin{theor}
	Soit $ f $ : $ \Hh^{2} \rightarrow \Hh^{2} $ une application 
	harmonique $ C^{2} $.
	\begin{itemize}
	    \item $ f $ est $ C^\infty $
	    \item $ f $ v\'erifie les estim\'ees de Schauder :    
	          $$ |\nabla^k f(x)| \leq C_{k,r} \sup B_{x,r} |df| $$
	\end{itemize}
    \end{theor}

    Lorsque $ f $ est d'image born\'ee, on a \'egalement des estim\'ees \`a 
    priori :
    
    \begin{theor}
	\label{principe-de-bloch}
	Soit $ f $ : $ \Hh^{2} \rightarrow \Hh^{2} $ une application 
	harmonique $ C^{2} $.
	Soit $ x \in \Hh^{2} $, $ r $ et $ R $ deux r\'eels positifs tels que 
	l'image de la boule $ B_{x,r} $ est incluse dans $ B_{f(x),R} $.
	Toutes les d\'eriv\'ees de $ f $ en $ x $ sont control\'ees par $ r $
	et $ R $ :
	$$ |\nabla^k f(x)| \leq K(r,R) $$
    \end{theor}
    
    Les propri\'et\'es de stabilit\'es des applications harmoniques du plan 
    hyperboliques se d\'eduisent de la formule de Bochner pour la 
    distance :
    
    \begin{theor}
	Soient $ f $ et $ g $ : $ \Hh^{2}\rightarrow\Hh^{2} $.
	La distance entre $ f $ et $ g $ v\'erifie l'in\'egalit\'e de Bochner :
	$$ \Delta d(f,g) \geq -|\tau(f)|-|\tau(g)|+\tanh(d(f,g)) $$
    \end{theor}

    La formule de Bochner pour la distance permet de d\'efinir des 
    crit\`eres d'existence d'applications harmoniques.
    
    \begin{deft}
	Soit $ f $ une application $ C^2 $ de $ \Hh^{2} $ dans $ \Hh^2 $.
	$ f $ est asymptotiquement quasi harmonique s'il existe un compact
	$ K $ du plan hyperbolique et un r\'eel $ \epsilon > 0 $ tel que, pour 
	tout $ x \in \Hh^{2}\setminus K $, $ f $ v\'erifie :
	$$ \lambda_{2}(f)^{2} \geq (1+\epsilon)|\tau(f)|+\epsilon $$
    \end{deft}
    
    D'apr\`es la formule de Bochner pour la distance, on a le th\'eor\`eme 
    (d\'ej\`a remarqu\'e par Li et Tam \cite[Theorem 6.4]{Li-Tam-93}) :
    
    \begin{theor}
	\label{existence-quasi-harmonique}
	Soit $ f $ une application asymptotiquement quasi harmonique du plan 
	hyperbolique. Il existe une application harmonique \`a distance born\'ee 
	de $ f $.
    \end{theor}
    
    \pv
    Fixons $ o\in \Hh^2 $ et $ R $ tel que pour tout $ x $ de
    $ \Hh^2 \setminus B_{o,R} $, on ait l'in\'egalit\'e :
    $$ \lambda_{2}(f)^{2} \geq (1+\epsilon)|\tau(f)|+\epsilon $$
    Pour tout $ r\in \R $, fixons $ K_{r} = f^{-1} B_{f(o),r} $.
    D'apr\`es le th\'eor\`eme de Schoen et Yau \cite{Schoen-Yau-78},
    il existe une unique application harmonique $ f_{r} $ de $ K_{r} $ 
    dans $ B_{f(o),r} $ qui co\"{\i}ncide avec $ f $ sur le bord de
    $ K_{r} $.
    
    Si $ \psi $ est une application $ C^2 $ et croissante
    de $ \R $ dans $ \R $, on peut d\'efinir :
    $ \phi(x)=d(f_r(x),f(x))+\psi(\rho(x)^2) $, o\`u $ \rho(x)=d(o,x) $.
    Le laplacien de $ \phi $ se calcule par~:
    \begin{eqnarray*}
	\Delta \phi
	& = & \Delta d(f_r,f)
	+ 2( |d\rho|^2+\rho\Delta \rho )\dot{\psi}(\rho) 
	+ 4\rho^2\ddot{\psi}(\rho)  \\
	& \geq  & \Delta d(f_r,f)
	+ 2\dot{\psi}(\rho)
	+ 4\rho^2\ddot{\psi}(\rho) 
    \end{eqnarray*}
    car on a fait l'hypoth\`ese que $ \psi $ est croissante
    et on sait que $ |d\rho|=1 $ et que $ \Delta\rho \geq 0 $
    puisque $ M $ est \`a courbure n\'egative.
    
    On sait que, en dehors de $ B_{o,R} $~:
    \begin{eqnarray*}
	\Delta d(f_r,f)
	& \geq  & \lambda_{2}(f)^2 \tanh(d(f_r,f))-|\tau(f)| \\
	& \geq  & (1+\epsilon)|\tau(f)|
	\left(\tanh(d(f_r,f))-\frac{1}{1+\epsilon}\right)
	+ \epsilon\tanh(d(f_r,f))
    \end{eqnarray*}
    
    On peut choisir $ \psi(x) $ tel que~:
    \begin{itemize}
	\item sur le compact $ B_{o,R} $~: $ \Delta \phi > 0 $.
	Si $ |\tau(f)|\leq C $ sur $ K $, il suffit en effet de choisir
	$ \psi(x) = (C+1)x $ pour $ x\leq R^2 $,
	\item en dehors de $ B_{o,R} $, lorsque $ d(f_r,f) \geq 1 $~:
	$$ \Delta \phi > (1+\epsilon)|\tau(f)|
	\left(\tanh(d(f_{r},f))-\frac{1}{1+\epsilon}\right) $$
	Si on pose $ \bar{\epsilon}=\epsilon \tanh(1)  $,
	il suffit en effet de choisir $ \psi $
	de sorte que $ \bar{\epsilon}+4x\ddot{\psi}(x) > 0 $
	et $ \dot{\psi}(x)\geq 0 $ lorsque $ x\geq R^2 $.
    \end{itemize}
    On a donc les contraintes suivantes :
    $ \dot{\psi}(R^2)=C+1>0 $ 
    et $ \ddot{\psi}(x) > -\frac{\bar{\epsilon}}{4x} $.
    Comme $ \int_t^\infty \frac{dx}{x} =\infty $,
    on peut choisir $ \psi $ de telle sorte que
    $ \psi(x) $ soit constante pour $ x\geq M(\epsilon,C) $.
    Dans ce cas, pour tout $ x\geq 0 $, on aura ~:
    $$ 0\leq \psi(x)\leq \alpha=\alpha(\epsilon,C) $$

    On sait que $ \phi $ atteint son maximum sur $ K_{r} $.
    Comme $ \Delta \phi > 0 $ sur $ B_{o,R} $,
    le maximum est atteint sur l'int\'erieur $ K_r\setminus B_{o,R} $, 
    ou sur le bord de $ K_{r} $.
    Trois cas sont donc \`a consid\'erer~:
    \begin{itemize}
	\item $ \phi $ atteint son maximum sur le bord de $ K_{r} $.
	On sait que $ \phi \leq d(f_r,f)+\sup(\psi) $ et $ \sup(\psi)=\alpha $,
	donc~: 
	$$ \phi \leq \alpha $$
	\item $ \phi $ atteint son maximum en un point $ x $
	int\'erieur \`a $ K_r $,
	et $$ d(f_r(x),f(x))\leq 1 $$
	Comme $ \phi\leq d(f_r,f)+\alpha $, on a~:
	$$ \phi\leq \alpha + 1 $$
	\item $ \phi $ atteint son maximum en un point $ x $
	int\'erieur \`a $ K_r $,
	et $ d(f_r(x),f(x))\geq 1 $.
	On sait que $ x\not\in B_{o,R} $, donc que :
	$$ \Delta \phi > (1+\epsilon)|\tau(f)|
	\left(\tanh(d(f_{r},f))-\frac{1}{1+\epsilon}\right) $$
	Comme $ \phi $ atteint son maximum en $ x $, on sait
	que $ \Delta \phi(x) \leq 0 $,
	et on en d\'eduit~:
	$$ d(f_r(x),f(x)) \leq \tanh^{-1}\left(\frac{1}{1+\epsilon}\right) $$
    \end{itemize}
    Dans tous les cas~: 
    $$ \max_{x\in K_{r}} \phi(x)
    \leq \tanh^{-1} \left(\frac{1}{1+\epsilon}\right) + \alpha + 1 $$
    Comme~: $ d(f_r,f) = \phi - \psi \leq \phi $, on en d\'eduit que
    $ f_r $ reste \`a distance born\'ee de $ f $~:
    $$ d(f_r,f) \leq \tanh^{-1} \left(\frac{1}{1+\epsilon}\right) + \alpha + 1 $$
    D'apr\`es le th\'eor\`eme \ref{principe-de-bloch},
    on sait donc que $ f_r $ est localement $ C^\infty $ \'equicontinue,
    born\'ee et on en conclut que $ f_r $ converge, \`a extraction pr\`es,
    vers $ f_\infty $ harmonique, \`a distance born\'ee de $ f $.
    \vp
    
    \begin{rem}
	Comme je l'ai explicit\'e dans ma th\`ese \cite{Rivet-99}, la d\'efinition 
	des applications asymptotiquement quasi harmonique et le th\'eor\`eme 
	d'existence d'applications harmoniques \`a distance born\'ees des 
	applications asymptotiquement quasi harmonique peut se 
	g\'en\'eraliser au cas o\`u $ M $ est simplement connexe, \`a courbure 
	n\'egative ou nulle et $ N $ v\'erifie une hypoth\`ese de courbure 
	strictement n\'egative.
    \end{rem}
    
    \section{M\'etriques harmoniques}
    
    Dans cette section, nous \'etudions les m\'etriques harmoniques 
    hyperboliques des surfaces de Riemann. La notion de m\'etrique 
    harmonique permet d'\'etudier les submersions harmoniques en 
    s'affranchissant de l'action du groupe des isom\'etries au but,
    en substituant \`a l'\'etude de $ f $ : $ M\rightarrow N $
    celle de la m\'etrique induite $ f^\ast g $. Lorsque $ f $ 
    est une submersion, la connaissance de la m\'etrique $ f^\ast g $ 
    est \'equivalente \`a la connaissance de $ f $, \`a l'action pr\`es des 
    isom\'etries de $ N $.
    
    \begin{deft}[M\'etrique harmonique].\\
	Soit $ (M,g_0) $ une surface de Riemann.
	Une m\'etrique $ g $ sur $ M $ est harmonique par rapport \`a
	$ g_0 $ si l'application identit\'e~:
	$ (M,g_0)\rightarrow(M,g) $ est harmonique.
    \end{deft}
    
    Si $ g $ est une m\'etrique sur $ (M,g_0) $, sa diff\'erentielle de 
    Hopf $ \phi $ d\'ecrit la partie sans trace de $ g $ par rapport \`a
    $ g_{0} $ et l'\'energie $ e $ la trace de $ g $ par rapport \`a $ g_0 $ :
    $$ g = \phi + e \times g_0 + \bar{\phi} $$
    Si on pose $ h = \log(\frac{e+\sqrt{e^2-4|\phi|^2}}{2}) $,
    l'\'energie peut se d\'ecomposer en $ e = e^h + e^{-h}|\phi|^2 $.
    La m\'etrique $ g $ est positive si et seulement si :
    $ e^2-4|\phi|^2 \geq 0 $. $ h $ est donc bien d\'efinie d\`es lors 
    que $ g $ est une m\'etrique d\'efinie positive et
    $ g $ peut alors s'\'ecrire :
    $$ g = \phi + (e^h +|\phi|^2 e^{-h})  g_0 + \bar{\phi} $$
    Dans toute la suite de cette section, nous supposerons les 
    m\'etriques $ g $ d\'efinies positives, ce qui revient \`a imposer :
    $$ h > \log|\phi| $$
    Sous cette hypoth\`ese, le th\'eor\`eme de caract\'erisation des 
    immersions harmoniques par leur diff\'erentielle de Hopf devient :
    
    \begin{theor}
	Soit $ (M,g_0) $ une surface de Riemann. Une m\'etrique $ g $ sur
	$ M $ est harmonique par rapport \`a $ g_0 $ si et seulement
	sa diff\'erentielle de Hopf est holomorphe.
    \end{theor}
    
    A une m\'etrique harmonique
    $ g $ sont naturellement associ\'es~:
    \begin{itemize}
	\item {\bf La diff\'erentielle de Hopf $ \phi $}.
	      Cette diff\'erentielle d\'efinit une m\'etrique plate
	      $ g_{|\phi|} = |\varphi||dz|^2 $,
	      et une paire de feuilletages,
	      avec des singularit\'es aux points o\`u $ \phi=0 $.
	      Le {\bf feuilletage horizontal} est caract\'eris\'e par 
	      $ \phi(X_h,X_h)\in \R^+ $,
	      le {\bf feuilletage vertical} par~: $ \phi(X_v,X_v)\in \R^- $.
	      Si $ \phi(x)\neq 0 $, il existe une carte dans laquelle 
	      $ \phi=dz^2 $.
	      Le feuilletage horizontal est d\'efini par les courbes 
	      $ y=C^{te} $,
	      le feuilletage vertical par les courbes $ x=C^{te} $.
	      Si on pose $ u = h-\log|\phi| $,
	      la m\'etrique harmonique s'\'ecrit~:
	      $$ g=4\left[\cosh^2(\frac{u}{2})dx^2
	      +\sinh^2(\frac{u}{2})dy^2\right] $$
	      Le feuilletage horizontal correspond aux directions de 
	      dilatation maximale~:
	      $ (\lambda_1)^2 = e^h(1+e^{-u})^2 $,
	      et le feuilletage vertical aux directions de dilatation 
	      minimale~:
	      $ (\lambda_2)^2 = e^h(1-e^{-u})^2 $,
	\item {\bf Le facteur de distortion quasi conforme 
	      $ u = h-\log|\phi| $}.
	      
	      $ u=\infty $ lorsque $ \phi=0 $, c'est-\`a-dire aux points o\`u la 
	      m\'etrique harmonique $ g_{h,\phi} $ est conforme;
	      et $ u=0 $ aux points o\`u $ g_{h,\phi} $ est de rang un.
	      
	      Il existe deux fonctions $ h_\pm $ telles que~: 
	      $ g=g_{h,\phi} $,
	      qui se d\'eduisent l'une de l'autre par~: 
	      $ h_+=2\log|\phi|-h_- $.
	      Les coefficients de distortion associ\'es sont $ u_+=-u_- $.
	      
	      Si $ g_{h,\phi} $ n'est pas d\'eg\'en\'er\'ee et $ \phi\neq 0 $,
	      $ u $ est d\'efini d\`es que $ \phi(z)\neq 0 $, donc sur un 
	      ouvert connexe dense et $ u $ ne s'annule jamais.
	      Quitte \`a remplacer $ h $ par $ 2\log|\phi|-h $ (donc $ u $ 
	      par $ -u $), on supposer dans ce cas que~: $ u > 0 $,
        \item {\bf La m\'etrique conforme $ g_h=e^h g_0 $}.
	      Les coefficients de dilatation de $ g_{h,\phi} $ sont~: 
	      $ \lambda_\pm^2 = e^h(1\pm e^{-u})^2 \leq 4 e^h $ ;
	      d'o\`u l'on d\'eduit~: $ g_{h,\phi}\leq 4 g_h $,
	      et $ g_h $ est (\`a une constante pr\`es) la plus petite 
	      m\'etrique conforme majorant $ g_{h,\phi} $,
	\item {\bf Le coefficient de Beltrami} de la m\'etrique est
	      $ \mu = e^{-u}\frac{\bar{\phi}}{|\phi|} $.
    \end{itemize}
    
    Ces diff\'erents objets g\'eom\'etriques sont au coeur des propri\'et\'es 
    des m\'e\-tri\-ques harmoniques hyperboliques :
    \begin{itemize}
	\item {\bf Unicit\'e et principe du maximum :}
	      il existe une unique m\'etrique harmonique hyperbolique $ g $
	      de diff\'erentielle de Hopf $ \phi $ telle que la m\'etrique
	      conforme associ\'ee est compl\`ete. De plus, si $ \dot{g} $ est une 
	      autre m\'etrique harmonique qui a la m\^eme diff\'erentielle de Hopf, 
	      on a l'in\'e\-ga\-li\-t\'e : $ \dot{g} < g $,
	\item {\bf Existence :} si $ \phi $ est une diff\'erentielle 
	      quadratique holomorphe sur $ (M,g_0) $, il existe une unique 
	      m\'etrique harmonique hyperbolique $ g $
	      de diff\'erentielle de Hopf $ \phi $ 
	      telle que la m\'etrique conforme associ\'ee est compl\`ete,
	\item {\bf Compl\'etude :} si $ \phi $ est une diff\'erentielle 
	      quadratique holomorphe bor\-n\'ee sur $ \Hh^2 $, la m\'etrique 
	      harmonique hyperbolique associ\'ee est compl\`ete et quasi 
	      conforme,
	\item {\bf D\'eg\'en\'erescence :} si $ \phi=\alpha dz^{2} $ sur le demi 
	      plan de Poincar\'e $ \U $, la m\'etrique harmonique hyperbolique 
	      associ\'ee est compl\`ete si et seulement si $ \alpha\in\R^- $. 
	      Lorsque $ \alpha\not\in \R^- $, les feuilles du feuilletage 
	      vertical ont une extr\'emit\'e de longueur finie, et les feuilles 
	      du feuilletage horizontal convergent le long de cette extr\'emit\'e
	      vers une g\'eod\'esique.
    \end{itemize}
    
    Wan a d\'emontr\'e dans \cite{Wan-92} le principe du maximum et les 
    th\'eor\`emes d'existence et d'unicit\'e pour les diff\'erentielles 
    quadratiques holomorphes born\'ees sur $ \Hh^2 $.
    Il a \'etendu les th\'eor\`emes 
    d'existence pour les diff\'erentielles quadratiques holomorphes non 
    born\'ees sur $ \Hh^{2} $ et sur $ \C $ dans \cite{Au-Wan}.
    Les articles de Han, Tam, Treibergs et Wan 
    \cite{Han-Tam-Treibergs-Wan}
    et de Shi, Tam et Wan \cite{Shi-Tam-Wan-99}
    donnent des exemples de m\'etriques harmoniques sur $ \C $ et 
    montrent que la g\'eom\'etrie de ces m\'etriques harmoniques est 
    intimement li\'ee aux feuilletages horizontaux et verticaux 
    associ\'es \`a la diff\'erentielle de Hopf.
    Notons que l'\'etude de la compl\'etude des m\'etriques harmoniques 
    hyperboliques reste encore mys\-t\'e\-rieu\-se : il existe de nombreux 
    exemples de m\'etriques hyperboliques compl\`etes dont la 
    dif\-f\'e\-ren\-tiel\-le de Hopf n'est pas born\'ee : Li et Tam 
    \cite{Li-Tam-93} ont donn\'e des exemples de dif\-f\'eo\-mor\-phis\-mes 
    harmoniques du plan hyperbolique qui ne sont pas quasi conformes. 
    D'apr\`es la th\'eorie de Wan, leur diff\'erentielle de Hopf n'est pas 
    born\'ee. La question de la compl\'etude des m\'etriques harmoniques 
    hyperboliques sur $ \C $ n'est pas tranch\'ee. Schoen a conjectur\'e 
    \cite{Schoen-93} qu'il n'existe pas de diff\'eomorphisme harmonique 
    de $ \C $ dans $ \Hh^{2} $, c'est \`a dire qu'il n'existe pas de 
    m\'etrique harmonique hyperbolique compl\`ete sur $ \C $.
    Cette conjecture a \'et\'e v\'erifi\'ee pour des cas particuliers :
    Han, Tam, Treibergs et Wan \cite{Han-Tam-Treibergs-Wan} ont 
    d\'emontr\'e que la conjecture de Schoen est v\'erifi\'ee lorsque la 
    diff\'erentielle de Hopf est polynomiale. Ce r\'esultat a \'et\'e \'etendu 
    pour d'autre classes de diff\'erentielles quadratiques holomorphes 
    par Au, Wan et Tam \cite{Au-Tam-Wan},
    mais la conjecture de Schoen reste encore ouverte dans 
    le cas g\'en\'eral.
    
    Avant de d\'etailler la th\'eorie de Wan, commen\c{c}ons par rappeler 
    l'expression de la courbure d'une m\'etrique harmonique :
    
    \begin{prop}
	Soit $ (M,g_0) $ une surface de Riemann,
	$ \phi $ une diff\'erentielle quadratique holomorphe
	et $ h $ une fonction $ C^\infty $ telle que
	$ h > \log|\phi| $.
	La m\'etrique $ g_{h,\phi} = \phi
                	+ ( e^h + |\phi|^2 e^{-h}) g_0
			+ \bar{\phi} $
	a pour courbure~:
	$$ K_{g_{h,\phi}} =
	                {-\frac{1}{2}\Delta h
			+ K_{g_0}}{e^h - |\phi|^2 e^{-h}}  $$
    \end{prop}

    \pv
    Supposons $ \phi\neq 0 $,
    et choisissons une carte conforme, dans laquelle $ \phi=dz^2 $.
    Le laplacien associ\'e \`a $ g_0(x)=\sigma^2(x) (dx^2+dy^2) $ est~:
    $$ \Delta = \frac{1}{\sigma}(\frac{\partial^2}
    {\partial x^2}+\frac{\partial^2}{\partial y^2}) $$
    Si $ u = h - \log|\phi| $ est le facteur de distortion quasi 
    conforme de la m\'etrique $ g_{h,\phi} $, on sait que :
    $$ g_{h,\phi}=4\left[\cosh^2(\frac{u}{2})dx^2
	         +\sinh^2(\frac{u}{2})dy^2\right] $$
    Posons $ \alpha = 4 \cosh^2(\frac{u}{2}) $ et
    $ \beta = 4 \sinh^2(\frac{u}{2}) $
    (en remarquant que $ d\alpha=d\beta=2\sinh(u)du $).
    En calculant explicitement les symboles de Christoffel, on 
    v\'erifie que la courbure de g est~:
    \begin{eqnarray*}
	K & = & - \frac{1}{2\alpha\beta}
	        \left( \alpha_{yy} + \beta_{xx}
		- \frac{1}{2}(\frac{\beta_x^2}{\beta} + \frac{\alpha_y^2}{\alpha}
		+ \frac{\alpha_x\beta_x}{\alpha} + \frac{\alpha_y\beta_y}{\beta}) 
		\right)          \\
	& = & - \frac{1}{2\alpha\beta}
	      \left( (\frac{\partial^2}{\partial 
	      x^2}+\frac{\partial^2}{\partial y^2})\alpha
	      - \frac{| d\alpha|^2 (\alpha+\beta)}{2\alpha\beta} 
	      \right)          \\
	& = & - \frac{1}{4\sinh(u)}(\frac{\partial^2}{\partial 
	x^2}+\frac{\partial^2}{\partial y^2})u                \\
	& = & - \frac{\Delta h+ \Delta \log \sigma}{2e^h - 2|\phi|^2 e^{-h}}
    \end{eqnarray*}
    Lorsque $ h=0 $, $ \phi=0 $, on retrouve par continuit\'e~:
    $ K_{g_0}=-\frac{1}{2}\Delta \log \sigma $
    et on en d\'eduit la formule annonc\'ee.
    \vp
    
    \subsection{Principe du maximum pour les m\'etriques har\-mo\-ni\-ques}

    Avant d'\'enoncer et de d\'emontrer le principe du maximum pour les 
    m\'e\-tri\-ques harmoniques, rappelons le principe du maximum d'Omori et 
    Yau :
    
    \begin{theor}
	Soit $ M $ une vari\'et\'e riemanienne compl\^ete, de courbure minor\'ee 
	et $ f $ une fonction d\'efinie sur $ M $.
	Si $ f $ est major\'ee, il existe une suite maximisante $ x_{n} $ telle 
	que :
	\begin{itemize}
	    \item la diff\'erentielle de $ f $ en $ x_{n} $ converge vers $ 0 $
	    $$ \lim_{n\rightarrow\infty}|df(x_{n}| = 0 $$
	    \item le laplacien de $ f $ en $ x_{n} $ est asymptotiquement n\'egatif 
	    ou nul :
	    $$ \lim_{n\rightarrow\infty} \Delta f(x_{n}) \leq 0 $$
	\end{itemize}
    \end{theor}
    
    Si $ g_{h,\phi} = \phi + (e^h+e^{-h}|\phi|^{2})g_{0}+\bar{\phi} $ 
    est une m\'etrique harmonique, $ g_{h,\phi} $ est hyperbolique si 
    et seulement si $ h $ est solution de l'\'equation elliptique :
    $$ \frac{1}{2}\Delta h = e^h-|\phi|^2 e^{-h}-1 $$
    et $ h > \log|\phi| $.
    En exploitant le principe du maximum d'Omori et Yau, on peut 
    v\'erifier que les m\'etriques harmoniques hyperboliques v\'erifient :
    
    \begin{theor}
	Soit $ \phi $ une diff\'erentielle quadratique holomorphe
	et $ g_{h,\phi} $ une m\'e\-trique harmonique hyperbolique
	telle que la m\'etrique conforme associ\'ee 
	$ g_h = e^h g $ est compl\`ete.
	Si $ g $ est une autre m\'etrique harmonique hyperbolique de m\^eme
	diff\'erentielle de Hopf $ \phi $,
	$ g < g_{h,\phi} $.
    \end{theor}
    
    % Cette propri\'et\'e est analogue au lemme de Schwarz pour les 
    % applications holomorphes
    
    \pv
    Si $ g $ a pour diff\'erentielle de Hopf $ \phi $, il existe une 
    fonction $ \ell $ : $ M\rightarrow \R $ telle que $ g $ s'\'ecrive 
    sous la forme $ g=g_{\ell,\phi} $ :
    $$ g_{\ell,\phi} = \phi
             	+ ( e^\ell + |\phi|^2 e^{-\ell}) g_0
		+ \bar{\phi} $$
    $ \ell $ et $ h $ sont solutions de l'\'equation elliptique :
    $$ \frac{1}{2}\Delta h = e^h-|\phi|^2 e^{-h}-1 $$
    Notons $ \Delta_h=e^{-h}\Delta $ le laplacien par rapport \`a
    la m\'etrique conforme $ g_h $.
    Si on pose $ \lambda = \ell - h $, on sait que~:
    \begin{eqnarray*}
	\frac{1}{2}\Delta_h \lambda & = & (e^\lambda-1)(1 + 
	e^{-2\bar{u}}e^{-2\lambda})
    \end{eqnarray*}
    Par hypoth\`ese, la m\'etrique conforme $ g_h $ est compl\`ete.
    De plus, la courbure de $ g_{h} $ est minor\'ee~: 
    $ K=-(1-e^{-2u})\geq -1 $.
    On peut donc appliquer le principe du maximum d'Omori et Yau sur
    $ (M,g_{h}) $.
    
    La fonction $ \alpha(x) = e^{-\frac{1}{3}\lambda(x)} $ est minor\'ee 
    par $ 0 $. D'apr\`es le principe du maximum d'Omori et Yau, il 
    existe une suite minimisante pour $ \alpha $ telle que :
    $$ \lim_{n\rightarrow\infty}d \alpha(x_n) = 0 $$ et 
    $$ \lim_{n\rightarrow\infty}\Delta_h\alpha(x_n) \geq 0 $$
    Si on suppose que $ \inf \alpha = 0 $, sachant que $ d\alpha(x_n) $ 
    tend vers $ 0 $, on en d\'eduit que
    \begin{eqnarray*}
	\Delta_h \alpha (x_{n})
	& \leq & - \frac{2}{3}e^{\frac{2}{3}w} + o(e^{\frac{2}{3}w}) \\
	& \rightarrow & - \infty
    \end{eqnarray*}
    ce qui contredit l'hypoth\`ese
    $ \lim_{n\rightarrow\infty}\Delta_h\alpha(x_n) \geq 0 $.
    Par la contrapos\'ee, on conclut que : $ \inf \alpha > 0 $, donc 
    que $ \lambda $ est major\'ee.
    
    On peut ainsi appliquer le principe du maximum \`a $ \lambda $ et en 
    d\'eduire qu'il existe une suite maximisante pour $ \lambda $ telle 
    que :
    $$ \lim_{n\rightarrow\infty}\Delta_h \lambda(x_n) \leq 0 $$
    Comme
    \begin{eqnarray*}
	\frac{1}{2}\Delta_\lambda w & = & (e^\lambda-1)(1 + 
	e^{-2\bar{u}}e^{-2\lambda})
    \end{eqnarray*}
    on en d\'eduit que~: $ e^{\sup \lambda}-1 \leq 0 $,
    c'est-\`a-dire~: $ \ell \leq h $.
    
    Si on pose $ f_{x}(y)=e^y+e^{-y}|\phi(x)|^2 $, on sait que :
    $$ g_{h,\phi}(x)-g_{\ell,\phi}(x) = [f_{x}(h(x))-f_{x}(\ell(x))] g_0 $$
    De plus :
    $ \frac{df_{x}}{dy} = e^y (1-e^{-v(y)}) $ o\`u $  v(y)=y-\log|\phi(x)| $.
    Si $ y \geq \log|\phi(x)| $, $ \frac{df_{x}}{dy} \geq 0 $.
    Comme $ h(x)\geq \ell(x) > \log|\phi(x)| $ pour tout $ x \in M $,
    on en d\'eduit que pour tout $ x \in M $, $ f_{x} $ est croissante 
    sur $ [\ell(x),h(x)] $. On en conclut que pour tout $ x \in M $,
    $ g_{\ell,\phi}(x) \leq g_{h,\phi}(x) $.
    \vp

    \subsection{M\'etriques harmoniques de diff\'erentielle de Hopf 
                born\'ee}
    
    Avant d'\'enoncer et de d\'emontrer le th\'eor\`eme d'existence de Wan 
    \cite{Wan-92}, nous rappelons le principe de prolongement 
    d'Aronszajn \cite{Aronszajn}, qui nous sera utile pour d\'emontrer 
    les propri\'et\'es de quasi isom\'etrie de la m\'etrique harmonique 
    d\'efinie par Wan :
    
    \begin{theor}
	Soit $ \CO $ un ouvert connexe de $ \Hh^2 $
	et $ u $~: $ \CO\rightarrow\R^+ $,
	tel que~: $ u(0)=0 $
	et $ 0\leq \Delta u \leq \kappa u $.
	Alors~: $ u = 0 $.
    \end{theor}
    
    \pv
    Si $  f(r)~:=\int_{B_r}u $,
    on a~: $  \dot{f}=\sinh(r)\varphi(r) $ ,
    o\`u~: $  \varphi(r)=\int_{|x|=1}u(rx) $ et 
    $ \varphi(0)=2\pi u(0) $ par continuit\'e.
    La formule de Green nous donne~:
    $ \sinh(r)\dot{\varphi} = \int_{B_r}\Delta u  $.
    On sait que~: $ 0\leq \Delta u \leq \kappa u $.
    Donc~: $ \dot{\varphi}\geq 0 $ et~:
    $  \sinh(r)\dot{\varphi} \leq \kappa\int_{B_r}u $.
    On en d\'eduit~:
    $  \dot{\varphi}  \leq  r\kappa\sup_{[0,r]}\varphi $.
    En particulier~:
    $ \dot{\varphi}(0)=0 $ par continuit\'e.
    Si $ u(0)=0 $, on a alors~:
    $ 0\leq \dot{\varphi} \leq r^2\kappa\sup_{[0,r]}\dot{\varphi} $.
    D'o\`u l'on conclut~: $ u = 0 $.
    \vp
    
    Le th\'eor\`eme d'existence de Wan s'\'enonce donc :
    
    \begin{theor}
	Soit $ \phi $ une diff\'erentielle quadratique holomorphe 
	bor\-n\'ee sur $ \Hh^2 $.
	Il existe une unique m\'etrique harmonique hyperbolique 
	compl\`ete $ g_\phi $,
	de diff\'eren\-tielle de Hopf $ \phi $.
	\begin{itemize}
	    \item $ g_\phi $ est maximale~: toute m\'etrique harmonique 
	          hyperbolique de m\^eme dif\-f\'e\-ren\-tiel\-le de
		  Hopf est major\'ee par $ g_\phi $,
	    \item $ g_\phi $ est quasi isom\'etrique \`a la m\'etrique de 
	          Poincar\'e,
	    \item $ g_\phi = \phi + (e^h+|\phi|^2 e^{-h})g_0 + \bar{\phi} $ 
	          o\`u $ h $ est l'unique solution born\'ee de~:
		  $$ \frac{1}{2}\Delta h = e^h-|\phi|^2 e^{-h}-1 $$
	\end{itemize}
    \end{theor}
    
    \begin{rem}
	Si $ \Omega $ est un domaine de $ \C $ contenant le disque
	$ \D=\{z~:|z|<1\} $ et si $ \phi $ est une diff\'erentielle 
	quadratique holomorphe born\'ee sur $ \Omega $,
	on peut associer
	\`a tout ouvert $ \Delta $ relativement compact dans $ \Omega $,
	une m\'etrique harmonique hyperbolique 
	$ g^\Delta_\phi $ de diff\'erentielle de Hopf  $ \phi $,
	qui est compl\^ete sur $ \Delta $ et quasi conforme.
	D'apr\`es le principe du maximum, cette famille 
	de m\'etriques
	est d\'ecroissante~: si $ \Delta\subset\Lambda $,
	$ g^\Lambda_\phi|_\Delta \leq g^\Delta_\phi $.
	
	Lorsque $ \phi=0 $, on retrouve la propri\'et\'e de monotonie
	pour la m\'etrique de Poincar\'e des domaines du plan 
	complexe.
	
	Cet exemple sugg\`ere d'interpr\'eter le principe du maximum
	pour les m\'e\-tri\-ques harmoniques hyperboliques
	comme une extension du lemme de Schwarz.
    \end{rem}

    \pv
    En appliquant la m\'ethode de Perron, on d\'emontre que l'\'equation :
    $$ \frac{1}{2}\Delta h = e^h-|\phi|^2 e^{-h}-1 $$
    admet une solution born\'ee.
    Posons $ F(h) = e^h - |\phi|^2 e^{-h} - 1 $
    et $ Q(h) = -\frac{1}{2}\Delta h + F(h) $.
    Si $ h^+=\log\frac{1+\sqrt{1+4\sup |\phi|^2}}{2} $ et $ h^-=0 $,
    on v\'erifie que $ Q(h^+)\geq 0 \geq Q(h^-) $.
    $ F $ \'etant croissante, on peut appliquer la m\'ethode de Perron 
    (cf Wan \cite{Wan-92} pour les d\'etails),
    et en d\'eduire qu'il existe une solution $ h\in[h^-,h^+] $.
    
    On sait que $ h $ est minor\'e et $ \log|\phi| $ est major\'e~:
    le facteur de distortion quasi conforme $ u=h-\log|\phi| $ est donc minor\'e.
    Quitte \`a faire agir le groupe des isom\'etries de $ \Hh^2 $, on 
    peut supposer que $ u $ atteint son minimum et que~:
    $ u(0)=\inf(u) $. En effet, si $ x_n $ est une suite minimisante 
    pour $ u $,
    posons $ (u_n,h_n,\phi_n)=(u,h,\phi)\circ\varphi_n $ o\`u
    $ \varphi_n $ est une isom\'etrie de $ \Hh^2 $ qui envoie $ 0 $ 
    sur le point $ x_n $.
    La suite $ (h_n,\phi_n)_{n\in\N} $ est $ C^\infty $ 
    \'equicontinue, donc converge (\`a extraction pr\`es)
    vers $ (h_\infty,\phi_\infty) $ tels que~: 
    $ u_\infty(0)=\inf (u_\infty) = \inf (u) $.
    $ u $ v\'erifie l'\'equation~:
    $$ \frac{1}{2}\Delta u = e^h(1-e^{-2u}) $$
    Comme $ u(0)=\inf(u) $,
    on a~: $ \Delta u(0)\geq 0 $,
    c'est-\`a-dire~: $ e^{h(0)}(1-e^{-2u(0)})\geq 0 $
    et donc $ u(0)\geq 0 $.
    
    Si on suppose que $ u(0)=0 $, on appliquer le principe de 
    prolongement d'Aronszajn sur un voisinage de $ 0 $. On sait en 
    effet qu'il existe un voisinage de $ 0 $ et $ \kappa > 0 $ tel que~:
    $$ 0\leq \frac{1}{2}\Delta u = e^h(1-e^{-2u})\leq \kappa u $$
    D'apr\`es le principe de prolongement d'Aronszajn, on en d\'eduit 
    que $ u $ s'annule sur un voisinage de z\'ero, puis, par connexit\'e, 
    que $ u $ s'annulle sur $ M $. En particulier, on en conclut que
    $ h = \log|\phi| $,
    et donc~: $$ \frac{1}{2}\Delta h = - 1 $$
    
    $ h $ \'etant minor\'ee, on peut appliquer
    le principe du maximum d'Omori et Yau
    et conclure qu'il existe $ x_n $ tel que~:
    $ \lim_{n\rightarrow\infty} \frac{1}{2}\Delta h(x_n) \geq 0 $,
    ce qui contredit l'hypoth\`ese $ \frac{1}{2}\Delta h = - 1 $.
    
    Par la contrapos\'ee, on en d\'eduit que $ \inf u > 0 $
    et que la m\'etrique $ g_{h,\phi} $ est quasi conforme.
    On sait de plus que les coefficients de dilatation de 
    $ g_{h,\phi} $ sont donn\'es par les formules :
    $ (\lambda_1)^2 = e^h(1+e^{-u})^2 $,
    et $ (\lambda_2)^2 = e^h(1-e^{-u})^2 $.
    $ h $ \'etant born\'e, on en d\'eduit que $ g_{h,\phi} $ est une quasi 
    isom\'etrie. De plus, la m\'etrique conforme associ\'ee $ g_h = e^h g_0 $,
    est compl\^ete. D'apr\`es le principe du maximum, toute m\'etrique 
    harmonique hyperbolique $ g $ de m\^eme diff\'erentielle de Hopf est 
    major\'ee par $ g_{h,\phi} $.
    \vp
    
    \vspace{1ex}
    
    La preuve du th\'eor\`eme d'existence nous donne des estim\'ees \`a 
    priori sur la fonction $ h $ :
    $$ 0 \leq h \leq \log\frac{1+\sqrt{1+4\sup |\phi|^2}}{2} $$
    En utilisant la majoration ainsi d\'emontr\'ee, Wan \cite{Wan-92}
    en d\'eduit une caract\'erisation des diff\'eomorphismes harmoniques 
    quasi conformes du plan hyperbolique :

    \begin{cor}
	Soit $ f $ un dif\-f\'eo\-mor\-phis\-me harmonique du plan 
	hyperbolique.
	$ f $ est quasi conforme si et seulement si $ f $ est une 
	quasi isom\'etrie.
    \end{cor}

    \pv
    Soit $ \phi $ la diff\'erentielle de Hopf de $ f $.
    Puisque $ f $ est un diff\'eomorphisme, la m\'etrique
    $ f^\ast g $ est une m\'etrique harmonique hyperbolique compl\`ete.
    La m\'etrique conforme associ\'ee $ g_h $ v\'erifie
    $ f^\ast g \leq 4 g_{h} $.
    En particulier, $ g_{h} $ est compl\`ete.
    
    Si $ \phi $ est born\'ee, d'apr\`es le th\'eor\`eme d'existence de Wan,
    $ f^\ast g $ est l'unique m\'etrique 
    harmonique hyperbolique maximale associ\'ee \`a $ \phi $.
    De plus, $ f^\ast g $ est une quasi isom\'etrie et la fonction $ h $ 
    associ\'ee v\'erifie :
    \begin{eqnarray*}
	h & \leq & h^+:=\log\frac{1+\sqrt{1+4\sup|\phi|^2}}{2} \\
	& \leq & \log(1+\sup|\phi|)
    \end{eqnarray*}
    Le facteur de distortion $ u =h-\log|\phi| $ est donc major\'e par :
    \begin{eqnarray*}
	u       & \leq & \log(1+\sup|\phi|)-\log|\phi|         \\
	\inf(u) & \leq & \log(1+\sup|\phi|)-\log\sup|\phi|     \\
	& \leq & \frac{1}{\sup|\phi|}
    \end{eqnarray*}
    
    Si $ \phi $ n'est pas born\'ee, on peut choisir une suites 
    exhaustive de compact de $ \Hh^{2} $ : par exemple
    $ \Omega_n=B_{x,n} $
    La diff\'erentielle quadratique holomorphe
    $ \phi $ est born\'ee pour la m\'etrique de
    Poincar\'e de $ \Omega_n $~:
    on peut lui associer une m\'etrique harmonique maximale $ g_n $.
    
    La m\'etrique $ f^\ast g $ restreinte \`a $ \Omega_n $ est une m\'etrique
    harmonique hyperbolique qui a m\^eme diff\'erentielle de Hopf que $ g_n $.
    D'apr\`es le principe du maximum~: 
    $ f^\ast g|_{\Omega_n}\leq g_n $.
    On en d\'eduit que les coefficients de distortions
    v\'erifient~: $ u\leq u_n $.
    
    En particulier :
    \begin{eqnarray*}
	\inf(u) & \leq &\inf(u_n)         \\
	& \leq & \frac{1}{\sup|\phi|_{\Omega_n}}
    \end{eqnarray*}
    
    Lorsque $ n\rightarrow\infty $, on en conclut~: $ \inf(u)=0 $,
    c'est \`a dire~: $ f $ n'est pas quasi conforme.
    \vp
    
    \subsection{Principe du maximum pour le facteur de distortion}
    
    Nous avons vu que les m\'etriques harmoniques v\'erifient un principe 
    du maximum, qui peut s'interpr\'eter comme une extension du lemme de 
    Schwarz. On peut \'egalement d\'emontrer que les facteurs de 
    distortion des m\'etriques harmoniques maximales v\'erifient un 
    principe du maximum :
    
    \begin{theor}
	Soient $ \phi_1 $ et $ \phi_2 $ deux diff\'erentielles quadratiques
	born\'ees du plan hyperbolique,
	$ g_1 $ et $ g_2 $ les m\'etriques harmoniques associ\'ees.
	Si $ |\phi_1(x)|\geq|\phi_2(x)| $ pour tout $ x\in\Hh^2 $,
	les coefficients de distortion $ u_1 $ et $ u_2 $
	des m\'etrique $ g_1 $ et $ g_2 $
	v\'erifient l'in\'egalit\'e inverse~:
	$$ u_1\leq u_2 $$
    \end{theor}

    Si on fixe une diff\'erentielle quadratique holomorphe $ \phi $, on
    en d\'eduit en particulier que les coefficients de distortion des 
    m\'etriques harmoniques hyperbolique associ\'ees \`a $ t\phi $ sont des 
    fonctions d\'ecroissantes de $ t $. On a ainsi une bonne image 
    g\'eom\'etrique de la d\'eg\'en\'erescence des m\'etriques $ g_{t\phi} $ :
    le feuilletage horizontal est dilat\'e, et le feuilletage vertical 
    est contract\'e lorsque $ |t| $ augmente.
    
    Dans le cas o\`u $ M $ est une surface de Riemann compacte, 
    Wolf \cite{Wolf-89} et Minsky \cite{Minsky-92} ont \'etudi\'e 
    diff\'erents cas de d\'eg\'en\'erescence des m\'etriques harmoniques. 
    Wolf \cite{Wolf-89} explicitait en particulier les liens entre la 
    compactification de Thurston de l'espace de Teichm\"uller des 
    surfaces compactes et les feuilletages horizontaux et verticaux 
    associ\'es aux diff\'erentielles de Hopf des m\'etriques harmoniques.
    
    % Remarquons que les principes du maximum pour les m\'etriques
    % et pour le facteur de distortion permettent de conclure que
    % l'on peut majorer localement les m\'etriques harmoniques 
    % hyperboliques et leur facteur de distortion.
    
    % Au contraire, trouver une minoration est un probl\`eme global.

    \vspace{1ex}

    \pv
    Posons 
    $ Q(h) = -\frac{1}{2}\Delta h + e^{h}-|\phi_2|^2 e^{-h} - 1 $
    et $ \lambda =u_1+|\phi_2| $.
    On v\'erifie que~:
    $$ Q(\lambda) = e^{h_2}(1-e^{-2u_1} ) 
    \left[ |\frac{\phi_1}{\phi_2}| - 1 \right] \geq 0 $$
    On sait donc que : $ 0 \leq \lambda $ et
    $ Q(0)\leq 0 \leq \lambda $
    En applicant la m\'ethode de Perron, on en d\'eduit qu'il existe
    $ \bar{h}_2 $ tel que : $ 0 \leq \bar{h_2} \leq \lambda $
    et $ Q(\bar{h_2})= 0 $.
    $ \bar{h_2} $ \'etant minor\'ee, la m\'etrique conforme
    $ e^{\bar{h_2}}g_0 $ est compl\`ete.
    
    D'apr\`es le principe du maximum pour les m\'etriques harmoniques 
    hyperboliques, on en d\'eduit~: 
    $ \bar{h_2} = h_2 $, puis $ h_2\leq \lambda=u_1+|\phi_2| $,
    c'est-\`a-dire~: $ u_2\leq u_1 $.
    \vp
    
    \subsection{M\'etriques harmoniques hyperboliques sur $ \Hh^2 $ et 
    sur $ \C $}
    
    Le th\'eor\`eme d'existence de Wan s'\'etend naturellement aux cas des 
    dif\-f\'e\-ren\-tiel\-les de Hopf non born\'ees sur $ \Hh^2 $ et aux cas des 
    diff\'erentielles de Hopf sur $ \C $. Dans sa plus grande 
    g\'en\'eralit\'e, le th\'eor\`eme est le suivant :
    
    \begin{theor}
	Soit $ (M,g) $ une surface de Riemann parabolique ou 
	hyperbolique et $ \phi $ une diff\'erentielle quadratique holomorphe.
	Il existe une unique m\'etrique harmonique hyperbolique 
	$$ g_\phi = \phi + (e^h+|\phi|^2 e^{-h})g_0 + \bar{\phi} $$
	telle que la m\'etrique conforme $ g_{h} = e^h g $ est compl\`ete.
	
	$ g_\phi $ est maximale~: toute m\'etrique harmonique 
	hyperbolique de m\^eme dif\-f\'e\-ren\-tiel\-le de
	Hopf est major\'ee par $ g_\phi $.
    \end{theor}
    
    Nous renvoyons aux articles de Wan et Au \cite{Au-Wan} dans le 
    cas parabolique et Wan et Tam \cite{Tam-Wan-94} dans le cas hyperbolique
    pour plus de d\'etails.
    
    Le principe du maximum pour le facteur de distortion reste 
    \'egalement valable dans ce cadre plus g\'en\'eral.

    La compl\'etude de $ g_{\phi} $ est un probl\`eme encore ouvert :
    \begin{itemize}
	\item Lorsque $ M=\Hh^2 $ et $ \phi $ n'est pas born\'ee, 
	      les deux cas : compl\'etude et non compl\'etude de la m\'etrique 
	      $ g_\phi $ peuvent se produire. Nous explicitons dans la section 
	      suivante un exemple frappant
	      (\'egalement \'etudi\'e par Shi, Wan et Tam \cite{Shi-Tam-Wan-99}).
	      Dans le mod\`ele du demi-plan, 
	      si $ \phi = \alpha dz^2 $ :
	      \subitem si $ \alpha \in \R^{-} $, la m\'etrique harmonique 
	               hyperbolique associ\'ee $ g_{\phi} $ est compl\`ete,
	      \subitem si $ \alpha \in \C\backslash\R^{-} $,
	               la m\'etrique harmonique hyperbolique associ\'ee
	               $ g_{\phi} $ n'est pas compl\`ete,
		       
	      Dans cet exemple, l'\'etude des feuilletages horizontaux et 
	      verticaux permet de d\'ecrire g\'eom\'etriquement les cas
	      o\`u $ g_{\phi} $ n'est pas compl\`ete,
	\item Lorsque $ M=\C $, Schoen a conjectur\'e que $ g_{\phi} $ 
	      n'est jamais une m\'etrique compl\`ete. Les travaux de Wan et al.
	      \cite{Tam-Wan-94} \cite{Shi-Tam-Wan-99}
	      ont permis de v\'erifier qu'il n'existe pas de contre-exemple 
	      simple \`a la conjecture de Schoen, en \'etudiant les cas o\`u 
	      $ \phi $ est un polyn\^ome, une exponentielle etc...
	      grace \`a l'\'etude des propri\'et\'es 
	      des feuilletages horizontaux et verticaux associ\'es \`a $ \phi $. 
	      Il n'existe pas \`a ma connaissance de preuve ni de contre 
	      exemple \`a la conjecture de Schoen.
    \end{itemize}
    
    \subsection{M\'etriques harmoniques invariantes sous l'action d'un 
    sous\-groupe \`a un param\`etre d'isom\'etrie}
    
    Si on cherche \`a classifier
    les diff\'erentielles quadratiques holomorphes 
    invariantes par un sous groupe \`a un param\`etre d'isom\'etries 
    puis les m\'etriques harmoniques associ\'ees,
    trois possibilit\'es sont \`a envisager~:
    \begin{itemize}
	\item $ \G $ est hyperbolique, donc conjugu\'ee \`a un sous groupe 
	      de la forme~:
	      $ z\rightarrow e^t z $ dans le mod\`ele du demi-plan.
	      Les diff\'erentielles invariantes par $ \G $ sont 
	      born\'ees~:
	      $ \phi=(\alpha+i\beta)\frac{dz^2}{z^2} $,
	\item $ \G $ est parabolique, donc conjugu\'ee \`a sous groupe de la 
	      forme~:
              $ z\rightarrow z+t $ dans le mod\`ele du demi-plan.
	      Les diff\'erentielles associ\'ees ne sont pas born\'ees~:
	      $ \phi=(\alpha+i\beta)dz^2 $.
	      Elles donnent des exemples int\'eressants
	      d'immersions harmoniques $ f $~: $ \Hh^2\rightarrow\Hh^2 $
	      (cf Li et Tam \cite{Li-Tam-93} Shi, Tam et Wan 
	      \cite{Shi-Tam-Wan-99}),
        \item $ \G $ admet un point fixe, donc est conjugu\'e \`a un 
	      sous groupe de la forme~:
	      $ z\rightarrow e^{it}z $ dans le mod\`ele du disque.
	      Les diff\'erentielles quadratiques holomorphes invariantes 
	      sous $ \G $
	      ont une singularit\'e au point fixe~: 
	      $ \phi=(\alpha+i\beta)\frac{dz^2}{z^2} $.
    \end{itemize}
    
    Nous calculons dans les sections suivantes les m\'etriques 
    harmoniques invariantes sous l'action d'un sous groupe 
    hyperbolique et parabolique. Les applications harmoniques associ\'ees
    ont \'et\'e calcul\'ee par Shi, Tam et Wan \cite{Shi-Tam-Wan-99} 
    par des m\'ethodes diff\'erentes de celle que nous exploitons ici.

    \subsubsection{Cas des sous groupes hyperboliques}
    
    Choisissons comme mod\`ele du plan hyperbolique
    l'espace $ \W $ des nombres de partie imaginaire comprises
    entre $ 0 $ et $ \pi $ 
    muni de la m\'etrique $$ g = \frac{1}{\sin^2(y)}(dx^2+dy^2) $$
    
    Les diff\'erentielles quadratiques holomorphe invariantes sous 
    l'action du sous groupe d'isom\'etrie
    hyperbolique $ z\rightarrow z+t $ sont de la forme~: 
    $ \phi=(\alpha+i\beta)dz^2 \in QD_b(\W) $ et le 
    dif\-f\'eo\-mor\-phis\-me harmonique associ\'e v\'erifie:
    $ f(z+t) = f(z) + \lambda t $.
    Il est donc de la forme~: 
    $ f(x+iy) = \lambda x + \varphi(y) + i\psi(y) $.
    
    La diff\'erentielle de Hopf $ \phi = (\alpha+i\beta)dz^{2} $
    d\'ecrit la partie sans trace de $ f^\ast g $.
    Par identification, on en d\'eduit~:
    \begin{eqnarray*}
	\lambda d\varphi & = & -2\beta\sin^2(\psi) dy \\
	d\psi            & = & \sqrt{\lambda^2-\frac{4\beta^2}
	         {\lambda^2}\sin^4(\psi)-4\alpha\sin^2(\psi)}dy
    \end{eqnarray*}
    $ f $ \'etant un diff\'eomorphisme, on sait que
    $ \int_0^\pi dy(\psi) = \pi $, c'est \`a dire~:
    $$ \int_0^\pi \frac{1}
    {\sqrt{\lambda^2-\frac{4\beta^2}
    {\lambda^2}\sin^4(\psi)-4\alpha\sin^2(\psi)}}
    d\psi = \pi $$
    $ \alpha $ et $ \beta $ \'etant fix\'es, on 
    peut calculer la constante $ \lambda(\alpha,\beta) $ par la 
    condition ainsi d\'efinie, puis les fonction $ \varphi $ et $ \psi $.
    
    On d\'ecrit ainsi tous les  dif\-f\'eo\-mor\-phis\-mes harmoniques
    quasi conformes, dont l'extension au bord du tube est une
    application quasi symm\'etrique, invariante sous l'action du 
    sous groupe hyperbolique
    $ z\rightarrow z+t $~:
    $ f(x)=\lambda x $ et $ f(x+i\pi)=(\lambda	x + i\pi) + K $.

    \subsubsection{Cas des sous groupes paraboliques}
    
    Le mod\`ele naturel pour \'etudier les diff\'erentielles quadratiques 
    ho\-lo\-mor\-phes invariantes sous l'action d'un sous groupe parabolique
    et les m\'etriques harmoniques hyperboliques associ\'ees
    est celui du demi plan de Poincar\'e :
    \begin{itemize}
	\item si $ \G $ est un sous groupe parabolique du groupe des 
	isom\'etries de $ \Hh^{2} $, il existe une isom\'etrie de $ \Hh^{2} $ 
	dans $ \U $ telle que l'action de $ \G $ sur $ \Hh^{2} $ est donn\'ee par
	$ z\rightarrow z+t $,
	\item dans cette carte, les diff\'erentielles quadratiques holomorphes 
	in\-va\-rian\-tes sous l'action de $ \G $ sont de la forme
	$ \phi=\alpha dz^{2} $
    \end{itemize}
    
    On peut calculer explicitement la m\'etrique harmonique hyperbolique 
    associ\'ee \`a $ \phi $ :
    
    \begin{theor}
	Dans le mod\`ele du demi plan $ \U=\{z|Im(z)>0\} $,
	si $ \phi=\alpha dz^2 $, la m\'etrique harmonique maximale associ\'ee \`a 
	$ \phi $ est~:
	$$ g_\phi =\frac{1+\cosh^2(2\sqrt{|\alpha|} y)}
	                {2\sinh^2(y)}(dx^2+dy^2)+\phi+\bar{\phi} $$
	\begin{itemize}
	    \item Si $ \alpha\in\R^- $, $ g_\phi $ est compl\`ete,
	    \item Si $ \alpha\not\in\R^- $, $ g_\phi $ n'est pas compl\`ete.
	          Les images des feuilles du feuilletage horizontal ont une 
		  extr\'emit\'e de longueur finie,
		  les points limites associ\'es d\'ecrivant une g\'eod\'esique.
	\end{itemize}
    \end{theor}

    Si $ \phi=\frac{1}{4}dz^2 $, on peut calculer l'immersion 
    harmonique associ\'ee \`a $ \phi $~:
    $$ f_\phi(x+iy)=e^x\left[\frac{1+i\sinh(y)}{\cosh(y)}\right] $$
    $ f_\phi $ n'est pas surjective~: l'image de $ f_\phi $ est le 
    quart de plan $ \{x+iy|x>0,y>0\} $.
   
    \vspace{1ex}

    \pv

    $ \phi $ est une diff\'erentielle quadratique non born\'ee,
    stable par le sous groupe d'isom\'etries paraboliques~: 
    $ z\rightarrow z+t $.
    La m\'etrique $ g_\phi $ associ\'ee est donc \'egalement stable et 
    l'immersion harmonique associ\'ee \`a $ \phi $ est invariante par un sous 
    groupe \`a un param\`etre~:
    $ f_\phi(z+t)=\gamma(t)f_\phi(z) $.
    
    Supposons que ce sous groupe est un sous groupe parabolique.
    Il est conjugu\'e \`a $ z\rightarrow z+\lambda t $,
    et $ f_\phi $ est de la forme~: 
    $ x+iy\rightarrow \lambda x+ g(y) $.
    Apr\`es action de l'homoth\'etie 
    $ z\rightarrow\lambda^{-1}z $), on peut r\'e\'ecrire
    $$ f_\phi(x+iy)=x+\rho(y)+i\varsigma(y) $$
    La diff\'erentielle de Hopf est~:
    \begin{eqnarray*}
	\phi & = & (\frac{1-\dot{\rho}^2-\dot{\varsigma}^2}{4\varsigma^2}
	- i\frac{\dot{\rho}}{2\varsigma^2})dz^2 \\
	& = & (a + ib) dz^{2}
    \end{eqnarray*}
    On en d\'eduit que~:
    \begin{eqnarray*}
	d\rho & = & 2\varsigma^2 b dy \\
	d\varsigma & = & \sqrt{1-4a\varsigma^2-4b^2\varsigma}dy
    \end{eqnarray*}
    Soit :
    \begin{eqnarray*}
	dy & = & \frac{d\varsigma}{\sqrt{1-4a\varsigma^2-4b^2\varsigma}}
    \end{eqnarray*}

    Comme $ \int dy(\varsigma)=\infty $, on conclut que~:
    $ a = 0 $, $ b = -\beta^{2} \leq 0 $, puis on calcule~:
    $ \rho(y)=0 $, 
    $ \varsigma(y)=\frac{1}{2\beta}\sinh(2\beta y) $.
    La m\'etrique $ g_\phi $ s'\'ecrit alors~:
    $$ g_\phi =\frac{1+\cosh^2(2\beta y)}
    {2\sinh^2(y)}(dx^2+dy^2)+\phi+\bar{\phi} $$
    et $ g_{\phi} = f^\ast g $ o\`u
    $$ f(x+iy) = x + i\frac{\sinh(2\beta y)}{2\beta} $$
    est un diff\'eomorphisme harmonique du plan hyperbolique.

    Si $ \phi = \alpha dz^{2} $ et $ \alpha\not\in \R^- $,
    on peut calculer explicitement la m\'etrique $ g_\phi $.
    En effet, $ g_{\phi} $ s'\'ecrit
    $ g_\phi=\phi+(e^h+|\phi|^2 e^{-h})g_0+\bar{\phi} $,
    o\`u $ h $ est l'unique solution de
    $ \frac{1}{2}\Delta h = e^h-|\phi|^2 e^{-h}-1 $
    telle que $ e^h g_0 $ est compl\`ete.
    On remarque que $ h $ ne d\'epend que de $ |\phi| $.
    D'apr\`es les calculs pour le cas $ \phi = -\beta^{2} dz^{2} $,
    on en d\'eduit la formule
    $$ g_\phi =\frac{1+\cosh^2(2\sqrt{|\alpha|} y)}
    {2\sinh^2(y)}(dx^2+dy^2)+\phi+\bar{\phi} $$
    
    Supposons que $ \phi=-\beta^2e^{2i\theta}dz^2 $.
    En posant $ u+iv=e^{-i\theta}(x+iy) $, la m\'etrique harmonique 
    hyperbolique associ\'ee \`a $ \phi $ est donc~:
	      
    \begin{eqnarray*}
	g_\phi & = & \frac{4\beta^2
	\left(du^2
	+ \cosh^2( 2\beta y) dv^2\right)}
	{\sinh^2 ( 2\beta y)}
    \end{eqnarray*}
	      
    Le feuilletage vertical est donn\'e par les courbes $ v=C^{te} $.
    L'extr\'emit\'e $ \gamma(t)=x(t)+iy(t) $ qui tend vers l'infini
    ( i.e~: $ y\rightarrow\infty $) a pour longueur~:
	      
    \begin{eqnarray*}
	\ell & = & \int_{u_0}^\infty \frac{2\beta dt}{\sinh(2\beta 
	t|\sin[\frac{\theta}{2}]|)}                 \\
	& \neq & +\infty
    \end{eqnarray*}

    La m\'etrique $ g_\phi $ n'est donc pas compl\`ete.
	      
    On sait de plus que le facteur de distortion $ u $ (et 
    toutes ses d\'eriv\'ees) tend vers z\'ero lorsque $ y\rightarrow\infty $.
    Dans une carte o\`u $ \phi=dz'^2 $,
    on peut v\'erifier que le feuilletage horizontal a pour courbure 
    g\'eod\'esique~:
    $$ K_h = -\frac{1}{4\cosh(\frac{u}{2})}\frac{du}{dy'} $$
    En particulier~: 
    $$ |K_h|\leq\frac{1}{4}\frac{|du|}{\sqrt{|\phi|}} $$
    Si $ (\Gamma_t,\gamma(t)) $ est la feuille horizontale point\'ee 
    passant par $ \gamma(t) $
    on en d\'eduit que la courbure de $ \Gamma_t $ converge 
    (uniform\'ement sur tout compact) vers $ 0 $,
    donc que $ \Gamma_t $ converge vers une g\'eod\'esique 
    $ \Gamma_\infty $.
    \vp

    \begin{rem}
	On voit sur cet exemple que la compl\'etude de la m\'etrique $ g_\phi $ 
	est li\'ee au comportement du feuilletage vertical de $ \phi $.
	On peut trouver de nombreux exemples de diff\'erentielles $ \phi $ 
	telles que la m\'etrique harmonique associ\'ee $ g_\phi $ n'est pas 
	compl\`ete, par exemple lorsque $ \phi=e^z dz^2 $ sur $ \U=\Hh^2 $
	ou, comme l'ont montr\'e Han Tam Treibergs et Wan 
	\cite{Han-Tam-Treibergs-Wan},
	lorsque $ \phi(z)=P(z)dz^2 $ si $ P $ est un polyn\^ome de degr\'e 
	$ d\geq 1 $ sur $ \U=\Hh^2 $.
	Dans ces deux cas, on peut trouver une feuille verticale qui diverge 
	dans la m\'etrique de Poincar\'e,
	mais qui est de longueur born\'ee pour la m\'etrique harmonique $ g_\phi $.
    \end{rem}
    
    \begin{rem}
	Li et Tam \cite{Li-Tam-93} ont construit des familles de 
	dif\-f\'eo\-mor\-phis\-mes harmoniques qui ne sont pas quasi conformes.
	Si on note fixe un point $ \alpha $ sur $ \partial_{\infty}\Hh^2 $, on lui 
	associe une famille $ f_{\alpha,t} $ de diff\'eomorphismes harmoniques.
	Dans le mod\`ele du demi plan, lorsque $ \alpha=\infty $, on a la 
	formule :
	$$ f_{\infty,t}(x+iy) = x+i\frac{\sinh(2t y)}{2t} $$
	Li et Tam ont d\'emontr\'e que si $ f $ est un diff\'eomorphisme 
	harmonique quasi conforme, $ \alpha_{1},\ldots,\alpha_{n} $ $ n $ 
	points sur $ \partial_{\infty}\Hh^{2} $,
	et $ \kappa_{1},\ldots,\kappa_{n} $ $ n $ nombres r\'eels,
	l'application
	$ f_{\alpha_n,\kappa_nt}\circ\ldots\circ f_{\alpha_1,\kappa_1t}\circ f $
	est aymptotiquement quasi-harmonique lorsque $ t $ est proche de $ 0 $.
	Elle est donc \`a distance born\'ee d'une application harmonique $ f_t $, qui 
	est une immersion lorsque $ t $ est suffisamment petit.

	A un point $ \alpha $ de $ \partial_{\infty}\Hh^{2} $
	est \'egalement associ\'ee une famille de 
	diff\'erentielle de Hopf $ \phi_{t}(\alpha_{n}) $
	telle que, dans le mod\`ele du demi plan,
	lorsque $ \alpha=\infty $, on a la 
	formule :
	$$ \phi_t(\infty) = - t^{2}dz^{2} $$
	Notons que cette d\'efinition est ambig\"ue,
	puisqu'il existe une famille \`a un 
	param\`etre d'isom\'etries de $ \Hh^{2} $ telle que
	$ g_{s}^\ast \phi_t(\infty) = \phi_{e^{2s}t}(\infty) $.
	
	On peut conjecturer que la diff\'erentielle de Hopf 
	des diff\'eomorphismes harmoniques d\'efinis par Li et Tam
	sont de la forme :
	$$ \phi = \phi_{0} + \sum_{k=1}^n \phi_{\lambda_{n}}(\alpha_{n}) $$
	o\`u $ \phi_{0} $ est une diff\'erentielle
	quadratique holomorphe born\'ee.
    \end{rem}

    \section{Extension des ho\-m\'eo\-mor\-phis\-mes qua\-si\-sym\-m\'e\-tri\-ques du cer\-cle}

    \subsection{Extension des $ C^1 $ diff\'eomorphismes du cercle}
    
    Nous rappelons dans ce paragraphe le th\'eor\`eme de Li et Tam
    \cite{Li-Tam-93}
    
    \begin{theor}[Li et Tam]
	Soit $ \partial f $ : $ \Sp^1\rightarrow\Sp^1 $ un $ C^1 $
	diff\'eomorphisme du cercle. $ \partial f $ s'\'etend sur $ \D=\Hh^2 $
	en $ f $, diff\'eomorphisme harmonique du plan hyperbolique. $ f $ et 
	toutes ses d\'eriv\'ees sont born\'ees :
	$$ \sup |df| \leq \varphi(\sup|d\partial f|) $$
    \end{theor}
    
    \begin{rem}
	On sait qu'il n'existe pas de m\'etrique riemanienne sur $ \Sp^{1} $
	invariante sous l'action du groupe d'isom\'etries de $ \Hh^{2} $.
	Si note $ \G $ le groupe d'isom\'etries de $ \Hh^{2} $,
	Li et Tam ont donc d\'emontr\'e un r\'esultat beaucoup plus fort :
	$$ \sup |df| \leq \varphi\left(
	\inf_{(\gamma_{1},\gamma_{2})\in \G^2}\left(
	\sup_{x\in\Sp^1}|d (\gamma_{1}\circ\partial f\circ\gamma_{2})(x)|
	\right)\right) $$
    \end{rem}
    
    \pv
    Fixons $ (x_{1},x_{2})\in \Sp^2 $. Il existe une partition de 
    l'unit\'e du disque ferm\'e $ (\chi_{1},\chi_{2}) $ telle que
    $ \chi_{i} = 1 $ sur un voisinage de $ x_{i} $.

    Pla\c{c}ons nous dans le mod\`ele du demi-plan : on identifie donc 
    l'espace de d\'epart $ (\D,x_{1},x_{2}) $ avec le demi plan 
    $ (\U, 0, \infty) $ et l'espace d'arriv\'ee
    $ (\D, \partial f(x_{1}), $ $ \partial f(x_{2})) $ avec le demi plan 
    $ (\U,0,\infty) $. 
    Dans cette carte, il existe une fonction
    $ h_{1} $ : $ \R \rightarrow \R $ $ C^1 $ \`a support compact qui co\"{\i}ncide 
    avec $ \partial f $ sur le support de $ \chi_{1} $.
    
    En utilisant les formules de Poisson, on sait que
    $ h_{1} $ s'\'etend en une fonction harmonique
    $ \varphi_{1} $ : $ \U \rightarrow \R $
    et $ |dh_{1}| $ s'\'etend en une 
    fonction harmonique $ \psi_{1} $.
    Notons
    $ f_{1}(x,y) = (\varphi_{1}(x,y) - y \psi_{1}(x);y\psi_{1}(x)) $.
    
    Li et Tam \cite{Li-Tam-93} ont d\'emontr\'e que $ f_{1} $ 
    a les propri\'et\'es suivantes :
    \begin{itemize}
	\item $ |\tau(f_{1})| $ tend vers $ 0 $ \`a l'infini,
	\item $ f_{1} $ est d'\'energie born\'ee et $ \lambda_{2}(f_{1}) $ reste 
	      loin de $ 0 $ au voisinage du support de $ \partial \chi_{1} $,
	\item si $ x $ appartient au support de $ \partial \chi_{1} $,
	      $ f_{1}(z) $ converge vers $ h_{1}(x) $ lorsque $ z $ tend 
	      vers $ x $.
    \end{itemize}
    
    On peut r\'ep\'eter la m\^eme construction sur le support de
    $ \chi_{2} $ et d\'efinir $ h(x) $ comme \'etant le barycentre des 
    points $ f_{i}(x) $ affect\'es des poids $ \chi_{i}(x) $.
    
    Li et Tam \cite{Li-Tam-93} d\'emontrent que l'application 
    ainsi d\'efinie est asymptotiquement quasi harmonique 
    et qu'elle s'\'etend en $ \partial f $ 
    sur le bord \`a l'infini de $ \Hh^2 $.
    D'apr\`es le th\'eor\`eme \ref{existence-quasi-harmonique},
    $ h $ est \`a distance born\'ee d'une 
    application harmonique $ f $.
    
    Dans le th\'eor\`eme \ref{existence-quasi-harmonique}, $ f $ est 
    limite d'une suite d'applications harmoniques
    $ f_{n} $ telle que $ f_{n} $ co\"{\i}ncide avec $ h $ sur le bord 
    d'une suite exhaustive de compact $ \Omega_{n}=h^{-1}(B_{x,n}) $.
    Lorsque $ n $ est suffisamment grand,
    $ f_{n}|_{\partial \Omega_{n}} = h_{n}|_{\partial \Omega_{n}} $
    est un hom\'eomorphisme \`a valeur dans le bord d'un convexe.
    D'apr\`es le th\'eor\`eme de Schoen et Yau \cite{Schoen-Yau-78}, $ f_n $ 
    est un diff\'eomorphisme.
    
    Lorsque $ n \rightarrow \infty $, on en d\'eduit que le coefficient 
    de distortion quasi conforme de $ f $ v\'erifie : $ u\geq 0 $.
    D'apr\`es le principe de prolongement d'Aronszajn, on a deux 
    possibilit\'es : soit $ u $ s'annule en un point $ x_{0} $, et dans 
    ce cas $ u=0 $ partout, soit $ u $ ne s'annule jamais.
    
    Dans le premier cas, on en d\'eduit que $ f $ est de rang $ 1 $ et 
    qu'il existe une application harmonique $ \alpha $ et une 
    g\'eod\'esique $ \gamma $ telles que $ f = \gamma \circ \alpha $.
    $ f $ est \`a distance born\'ee de $ h $ : $ f $ s'\'etend donc sur
    $ \partial_{\infty} \Hh^{2} $ et son extension sur le bord \`a 
    l'infini est l'hom\'eomorphisme $ \partial f $.
    Si $ f $ est de la forme $ \gamma \circ \alpha $, son extension 
    sur le bord \`a l'infini \`a une image r\'eduite aux deux extr\'emit\'es de 
    la g\'eod\'esique $ \gamma $.
    
    Par la contrapos\'ee, on en d\'eduit que $ u > 0 $ et que $ f $ est 
    une dif\-f\'eomor\-phis\-me.

    Par construction, $ h $ est d'\'energie born\'ee.
    D'apr\'es le principe de Bloch, on en d\'eduit que
    $ f $ est \'egalement d'\'energie born\'ee.
    La m\'etrique $ f^\ast g $ est donc une m\'etrique de la forme :
    $ f^\ast g = g_{h,\phi} $ o\`u $ \phi $ est une diff\'erentielle 
    quadratique holomorphe born\'ee.
    $ f^\ast g $ \'etant compl\^ete, la m\'etrique conforme associ\'ee l'est 
    \'egalement : $ f^\ast g $ est donc l'unique m\'etrique harmonique 
    hyperbolique maximale associ\'ee \`a $ \phi $.
    $ \phi $ \'etant born\'ee, $ f $ est une quasi isom\'etrie.
    \vp
    
    \subsection{Lemme de r\'egularit\'e}
    
    Dans le mod\`ele du disque de Poincar\'e, le bord \`a l'infini de 
    $ \Hh^2 $ s'identifie au cercle $ \Sp^1 $ de centre $ 0 $ et de 
    rayon $ 1 $. Une fois fix\'es trois points $ (x_{1};x_{2};x_{3}) $ sur 
    $ \partial_{\infty}\Hh^2 $, on peut identifier
    $ (\partial_{\infty}\Hh^2;x_{1};x_{2};x_{3}) $ au cercle
    $ (\Sp^1;1;e^{i\frac{2\pi}{3}};e^{-i\frac{2\pi}{3}}) $.
    $ (\partial_{\infty}\Hh^2;x_{1};x_{2};x_{3}) $ est ainsi
    naturellement muni d'une structure riemanienne. Sous ces 
    hypoth\`eses, on peut donc \'enoncer :
    
    \begin{lem}
	Soit $ \partial f $ : $ \Sp^1\rightarrow\Sp^1 $ un $ C^1 $ 
	diff\'eomorphisme $ K $ quasi symm\'etrique qui fixe les points 
	$ (1;e^{i\frac{2\pi}{3}};e^{-i\frac{2\pi}{3}}) $.
	Il existe une fonction $ \lambda $ : 
	$ \R^+\rightarrow\R^+ $ telle que
	la diff\'erentielle de $ f $ v\'erifie :
	$$ \frac{1}{\lambda(K)} \leq |df| \leq \lambda(K) $$
    \end{lem}
    
    \pv
    Compte tenu des sym\'etries du probl\`eme, il suffit de d\'emontrer :
    $ |df(e^{i\theta})| \leq \lambda(K) $ pour 
    $ \theta\in [0,2\frac{\pi}{3}] $. L'application 
    $$ g(z) = 
    -e^{i\frac{2\pi}{3}}\frac{z-1}{z-e^{-i\frac{2\pi}{3}}} $$
    est une isom\'etrie du disque de Poincar\'e vers le demi-plan qui 
    envoie les points $ (1;e^{i\frac{2\pi}{3}};e^{-i\frac{2\pi}{3}}) $ 
    sur $ (0;1;\infty) $ et qui est de d\'eriv\'ee born\'ee sur tout 
    compact de $ \D - \{e^{-i\frac{2\pi}{3}}\} $.
    
    $ f' = g\circ f\circ g^{-1} $ est un $ C^1 $ diff\'eomorphisme de $ 
    \R $ qui fixe $ (0,1) $.
    L'in\'egalit\'e \`a d\'emontrer pour $ f $
    est \'equivalente \`a l'in\'egalit\'e 
    pour $ f' $ :
    $ |df'(x)| \leq \lambda'(K) $
    pour tout $ x\in[0,1] $.
    
    $ f $ \'etant $ K $ quasi symm\'etrique, on sait que pour tout 
    quadrupl\'e $ (a,b,c,d) $ dans $ \R^4 $ on a :
    $$ \frac{f'(b) - f'(a)}{b-a}\times\frac{f'(c) - f'(d)}{c-d}\times
    \frac{c-a}{f'(c) - f'(a)}\times\frac{b-d}{f'(b) - f'(d)} \leq K $$
    Si on choisit $ a \in [0,1] $, $ b = a + t $, $ (c,d) = (-1;2) $, 
    on en d\'eduit, lorsque $ t\rightarrow 0 $ :
    $$ \frac{df'}{dt}(a) \leq K \frac{3}{f'(2) - f'(-1)}\times
    \frac{f'(a) - f'(-1)}{a+1}\times\frac{f'(2) - f'(a)}{2-a} $$
    En utilisant la quasi symm\'etrie de $ f $, on peut
    \'egalement d\'emontrer que
    $ 1+\frac{1}{K}\leq f'(2) \leq 2K $ et
    $ 2K+1 \leq f'(-1) \leq -\frac{1}{K} $.
    Sachant que $ 0\leq f'(a)\leq 1 $,
    on en d\'eduit qu'il existe bien une 
    fonction $ \lambda' $ telle que
    $ |df'(x)| \leq \lambda'(K) $.
    \vp
    
    \begin{cor}
	Soit $ \partial f $ un $ C^1 $ diff\'eomorphisme
	quasi symm\'etrique du bord de $ \Hh^2 $ et $ f $ l'extension 
	harmonique de $ \partial f $ sur $ \Hh^2 $.
	Il existe une fonction $ \varphi $ : $ \R^+\rightarrow \R^+ $
	telle que :
	si $ \partial_{\infty} f $ est $ K $ quasi symm\'etrique,
	$ f $ est une $ \varphi(K) $ quasi isom\'etrie.
    \end{cor}
    
    \pv
    Fixons $ (x_{1},x_{2},x_{3}) $ trois points sur
    $ \partial_{\infty}\Hh^2 $.
    On sait qu'il existe des isom\'etries $ g_{1},g_{2} $ :
    $ \Hh^2\rightarrow \D $ telles que
    $ g_{2}(x_{1},x_{2},x_{3})= 
    (1;e^{i\frac{2\pi}{3}};e^{-i\frac{2\pi}{3}}) $ et
    $ \partial f' = g_{1}\circ \partial f \circ g_{2} $ est un $ C^1 $ 
    diff\'eomorphisme de $ \Sp^1 $ qui fixe le tripl\'e
    $ (1;e^{i\frac{2\pi}{3}};e^{-i\frac{2\pi}{3}}) $.
    Si $ \partial f $ est $ K $ quasi symm\'etrique, d'apr\`es le lemme 
    de r\'e\-gu\-la\-ri\-t\'e, on d\'eduit que la diff\'erentielle de $ \partial f' $ 
    est born\'ee.
    D'apr\`es les estim\'ees $ C^1 $ de Li et Tam, on sait donc que 
    l'extension harmonique $ f $ de $ \partial f' $ a ses d\'eriv\'ees 
    born\'ees. On a donc une in\'egalit\'e :
    $$ \sup |df| \leq \varsigma_{0}\sup |d\partial f| $$
    D'apr\`es le th\'eor\`eme de Wan, $ f $ \'etant un diff\'eomorphisme 
    d'\'energie born\'ee, c'est une quasi isom\'etrie. Si on note $ K_{f} $ 
    le coefficient de quasi isom\'etrie de $ f $, on a une in\'egalit\'e :
    $$ K_{f} \leq \varsigma_{1} (\sup |df|) $$
    d'o\`u l'on d\'eduit l'in\'egalit\'e annonc\'ee
    en posant $ \varphi = \varsigma_1\circ\varsigma_0 $
    \vp
    
    \subsection{Extension des ho\-m\'eo\-mor\-phis\-mes qua\-si 
             sym\-m\'e\-tri\-ques}
    
    En exploitant le th\'eor\`eme d'extension des $ C^1 $ diff\'eomorphismes 
    du cercle et le lemme de r\'egularit\'e, nous pouvons d\'emontrer la 
    conjecture de Schoen :
    
    \begin{theor}
	Tout ho\-m\'eo\-mor\-phis\-me quasi sym\-m\'e\-tri\-que du 
	cer\-cle s'\'etend en un 
	dif\-f\'eo\-mor\-phis\-me harmonique du plan hyperbolique.
    \end{theor}
    
    \pv
    Soit $ \partial f $ un hom\'eomorphisme quasi symm\'etrique du 
    bord \`a l'infini de $ \Hh^2 $. D'apr\`es le th\'eor\`eme de Douady-Earle 
    \cite{Douady-Earle}, $ \partial f $ s'\'etend en un $ C^\infty $
    dif\-f\'eo\-mor\-phis\-me $ f_{DH} $ sur $ \Hh^2 $
    qui est une quasi isom\'etrie.
    
    Fixons $ O\in \Hh^2 $ et consid\'erons $ B_{r} $, la boule de 
    centre $ O $, de rayon $ r $, et $ C_{r} = f_{DH}^{-1}(B_{r}) $.
    $ B_{r} $ et $ C_{r} $ admettent chacun une m\'etrique de Poincar\'e
    $ g_{B_{r}} $ et $ g_{C_{r}} $.
    L'extension de Douady-Earle $ f_{DH} $ est $ K $ quasi conforme :
    elle s'\'etend donc en un hom\'eomorphisme $ \varphi(K) $ 
    quasi symm\'etrique de $ \partial_{\infty}(C_{r},g_{C_{r}}) $ vers
    $ \partial_{\infty} (B_{r},g_{B_{r}}) $.
    Notons que $ \varphi(K) $ ne d\'epend pas de $ r $.
    Par construction de $ B_{r} $ et $ C_{r} $,
    cet hom\'eomorphisme $ \partial_{\infty} f_{r} $ est $ C^\infty $.
    D'apr\`es le th\'eor\`eme d'extension $ C^1 $ de Li et Tam, il s'\'etend 
    en un diff\'eomorphisme harmonique $ f_{r} $.
    
    On d\'efinit ainsi une suite de fonctions $ f_{r} $ qui v\'erifient 
    les propri\'et\'es suivantes :
    \begin{itemize}
	\item $ f_{r} $ : $ (C_{r},g_{C_{r}})\rightarrow (B_{r},g_{B_{r}}) $
	      est une $ \psi(K) $ quasi isom\'etrie,
	\item $ f_{r} $ est \`a distance born\'ee de $ f_{DH} $ dans
	      $ (B_{r},g_{B_{r}}) $ :
	      $$ d(f_{r},f_{DH}) \leq \lambda(K) $$
    \end{itemize}
    o\`u l'on note que $ \psi(K ) $ et $ \lambda(K) $ ne d\'ependent pas 
    de $ r $. La premi\`ere propri\'et\'e se d\'eduit du lemme de r\'egularit\'e
    et la seconde est une propri\'et\'e classique de $ \Hh^2 $ :
    \begin{quotation}
	{\it 
	Si $ f_{1} $ et $ f_{2} $ sont deux hom\'eomorphisme $ K' $ quasi conforme
	de $ \Hh^2 $, $ f_{1} $ et $ f_{2} $ sont \`a distance born\'ee l'un de 
	l'autre, et la distance entre $ f_{1} $ et $ f_{2} $ est born\'ee 
	par $ \mu(K') $.}
    \end{quotation}

    On sait de plus que les m\'etriques $ g_{B_{r}} $ et $ g_{C_{r}} $
    sont d\'ecroissantes et con\-ver\-gent vers les m\'etriques de Poincar\'e 
    sur $ \Hh^2 $ et $ f_{DH}(\Hh^2) $.
    On en d\'eduit que la suite $ f_{r} $ est $ C^\infty $ \'equicontinue 
    sur tout compact de $ \Hh^2 $ et converge, \`a extraction pr\`es,
    pour la topologie $ C^2 $ sur tout compact, 
    vers $ f $, telle que :
    \begin{itemize}
	\item $ f $ : $ \Hh^2\rightarrow \Hh^2 $
	      est une $ \phi(K) $ quasi isom\'etrie,
	\item $ f $ est \`a distance born\'ee de $ f_{DH} $ dans
	      $ \Hh^2 $ :
	      $$ d(f_{r},f_{DH}) \leq \lambda(K) $$
	\item $ f $ est harmonique.
    \end{itemize}
    \vp
    
    \begin{cor}
	L'application de Wan qui associe
	\`a une diff\'erentielle quadratique holomorphe 
	$ \phi $ born\'ee sur $ \Hh^{2} $
	la classe de Teichm\"uller de la m\'etrique
	$$ g_{\phi} = \phi + (e^h+|\phi|^{2}e^{-h}) g_{\Hh^{2}} + \bar{\phi} $$
	o\`u $ h $ est l'unique solution born\'ee de :
	$$ \frac{1}{2}\Delta h = e^h-|\phi|^2 e^{-h}-1 $$
	d\'efinit une bijection de l'espace des diff\'erentielles quadratiques 
	holomorphes born\'ees dans l'espace universel de Teichm\"uller.
    \end{cor}
    
    Si on prend comme mod\`ele de l'espace universel de Teichm\"uller
    l'espace des hom\'eomorphismes quasi conformes de $ \Hh^{2} $,
    l'application de Wan associe \`a une diff\'erentielle quadratique 
    holomorphe born\'ee $ \phi $ la classe de Teichm\"uller des 
    diff\'eomorphismes harmoniques de diff\'erentielle de Hopf $ \phi $.
    
    Dans le mod\`ele o\`u l'espace universel de Teichm\"uller est l'espace 
    des hom\'eomorphismes quasi symm\'etriques du bord \`a l'infini de
    $ \Hh^2 $, l'inverse de l'application de Wan associe \`a un 
    hom\'eomorphisme quasi symm\'etrique $ \partial f $
    la diff\'erentielle de Hopf de l'extension harmonique quasi conforme 
    de $ \partial f $ sur $ \Hh^2 $.
    
    \bibliographystyle{plain}

\end{document}